\documentclass{article}

\usepackage{arxiv}

\usepackage[utf8]{inputenc} 
\usepackage[T1]{fontenc}    
\usepackage[hidelinks]{hyperref}       
\usepackage{url}            
\usepackage{booktabs}       
\usepackage{amsfonts}       
\usepackage{nicefrac}       
\usepackage{microtype}      
\usepackage{lipsum}         
\usepackage{graphicx}
\usepackage{doi}

\usepackage{colortbl}
\usepackage{listings}
\usepackage{bm}

\usepackage{amsmath}
\usepackage{cleveref}       

\usepackage[ruled,vlined]{algorithm2e}

\usepackage{appendix}
\usepackage{cite}

\usepackage{multirow}
\usepackage{array}
\newcommand{\head}[2]{\multicolumn{1}{>{\arraybackslash}p{#1}}{{#2}}}

\usepackage{chngcntr}

\SetKwComment{Comment}{/* }{ */}

\title{About optimal loss function for training\\ physics-informed neural networks\\ under respecting causality}

\date{\today}

\author{ \href{https://orcid.org/0000-0002-4930-1846}{\includegraphics[scale=0.06]{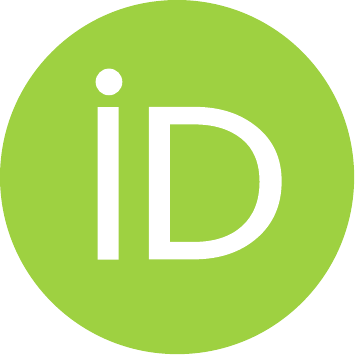}\hspace{1mm}Vasiliy A. Es'kin}\thanks{Corresponding author: Vasiliy Alekseevich Es’kin (vasiliy.eskin@gmail.com)} \\
	Department of Radiophysics, University of Nizhny Novgorod\\
	 Nizhny Novgorod, Russia, 603950\\
	 and\\
	 Manpower IT Solutions, Nizhny Novgorod, 603140, Russia\\
	\texttt{vasiliy.eskin@gmail.com} \\
	\And
	{\hspace{1mm}Danil V. Davydov} \\
	Mechanical Engineering Research Institute \\ 
	Russian Academy of Sciences\\
	Nizhny Novgorod, Russia, 603155\\	
	and\\
	Manpower IT Solutions\\
	Nizhny Novgorod, 603140, Russia\\
	\texttt{davidovdan27@yandex.ru} \\
	\And
	\href{https://orcid.org/0000-0001-9131-4256}{\includegraphics[scale=0.06]{orcid.pdf}\hspace{1mm}Ekaterina D. Egorova} \\
	Huawei Nizhny Novgorod Research Center\\
	Nizhny Novgorod, Russia\\
	and\\
	Institute of Applied Physics\\
	Russian Academy of Sciences\\
	Nizhny Novgorod, Russia, 603155\\	\texttt{egorovaed@yandex.ru} \\
	\And
	\href{https://orcid.org/0000-0001-5210-8281}{\includegraphics[scale=0.06]{orcid.pdf}\hspace{1mm}Alexey O. Malkhanov} \\
	Huawei Nizhny Novgorod Research Center\\
	Nizhny Novgorod, Russia\\	\texttt{alexey.malkhanov@gmail.com} \\
	\And
	{\hspace{1mm} Mikhail A. Akhukov } \\
	Manpower IT Solutions\\
	Nizhny Novgorod, 603140, Russia\\
	\texttt{ma.akhukov@yandex.ru} \\
	\And
	\href{https://orcid.org/0000-0002-0454-5249}{\includegraphics[scale=0.06]{orcid.pdf}\hspace{1mm}Mikhail E. Smorkalov} \\
	Skolkovo Institute of Science and Technology\\
	Moscow, Russia\\
	and\\
	Huawei Nizhny Novgorod Research Center\\
	Nizhny Novgorod, Russia\\	\texttt{smorkalovme@gmail.com} \\
}

\renewcommand{\headeright}{}
\renewcommand{\undertitle}{}
\renewcommand{\shorttitle}{A preprint}

\renewcommand{\vec}{\bf}

\hypersetup{
pdftitle={A template for the arxiv style},
pdfsubject={math.na, cs.AI, physics.comp-ph},
pdfauthor={Vasiliy A. Es'kin, Danil V. Davydov, Ekaterina D. Egorova, Alexey O. Malkhanov, Mikhail A. Akhukov, Mikhail E. Smorkalov},
pdfkeywords={Deep Learning, Physics-informed Neural Networks, Partial differential equations, Predictive modeling, Computational physics, Nonlinear dynamics},
}

\begin{document}
\maketitle

\begin{abstract}
	A method is presented that allows to reduce a problem described by differential equations with initial and boundary conditions to the problem described only by differential equations. The advantage of using the modified problem for physics-informed neural networks (PINNs) methodology is that it becomes possible to represent the loss function in the form of a single term associated with differential equations, thus eliminating the need to tune the scaling coefficients for the terms related to boundary and initial conditions. The weighted loss functions respecting causality were modified and new weighted loss functions based on generalized functions are derived. Numerical experiments have been carried out for a number of problems, demonstrating the accuracy of the proposed methods.
\end{abstract}


\section{Introduction}
The last decade has been notable by significant achievements in the field of machine learning such as deep learning~\cite{Alzubaidi2021}. The development of training methods and the improvement of deep neural network architectures have led to the creation of applications that are comparable or superior to the average person in many fields of activity previously available only to humans~\cite{Lecun2015,Linardatos2021,Khan2022}. Among these applications, machine translation systems~\cite{Vaswani2018,Wang2019,Yao2020}, image generators by description~\cite{Ramesh2021,FROLOV2021187}, game systems~\cite{Silver2018,Schrittwieser2020} and text generators~\cite{Ouyang2022} can be distinguished. The current revolution in the field of deep learning has not left aside the areas of human activities related to the obtaining and use of scientific knowledge. Auxiliary systems for solving scientific problems can be divided into those that are trained using previously obtained data (from experiments or computations) and build their predictions on this basis~\cite{Lu_2021,Li2022,Krinitskiy2022,Fanaskov2022,Ovadia2023}, those that leverage previously established physical laws~\cite{Raissi2019,Wang2022,Wang2022_2} and those that utilize the both approaches~\cite{Raissi2017, Raissi2017_2}. The first category includes the widely known Alphafold~\cite{Ronneberger2021} and, somewhat less famous, analysis system of interaction of laser radiation with the medium~\cite{Lin2022,Streeter2022}. Data-driven methods are out of consideration of this paper, and we rather focus on the second category exclusively.
Initial ideas for constraining neural networks using physical laws have been explored in \cite{Psichogios1992AHN,Lagaris_1998} studies of previous century. Though more contemporary research in the field was done few years ago in physics-informed neural networks (PINNs)~\cite{Raissi2019} with help of modern computational tools. In the PINN approach, a neural network is trained to approximate the dependences of physical values on spatial and temporal coordinates for a given problem described by physical equations together with initial and boundary conditions, and sometimes with data from computations or experiments. This method has been used in recent years to solve a wide range of problems described by ordinary differential equations~\cite{Raissi2017,Raissi2017_2,Raissi2019, Rackauckas2020}, integro-differential equations~\cite{Yuan2022}, nonlinear partial differential equations~\cite{Raissi2017,Raissi2017_2,Raissi2019,Jin2021,Zhao2021,Kharazmi2021}, PDE with noisy data~\cite{Yang2021}, etc.~\cite{Cuomo2022,Pang2020}, related to various fields of knowledge such as thermodynamics~\cite{Patel2022,Kharazmi2021}, hydrodynamics~\cite{Cai2021,Thakur2023,Jin2021,Zhao2021}, mechanics~\cite{Moseley2020}, electrodynamics~\cite{Lin2021,Kharazmi2021}, geophysics~\cite{He2020}, finance~\cite{Yuan2022}, etc.

Despite of the PINN success, it is often necessary to adapt hyperparameters of learning and the architecture of the neural network for each new problem. In particular, training a neural network boils down to minimization of the loss function, which includes terms related to differential equations, initial and boundary conditions~\cite{Raissi2017,Raissi2019,Wang2022}. In practice, for the training procedure to converge, it is necessary to tune the weighting coefficients for each of the terms mentioned. This tuning is done either empirically or by finding the optimum~\cite{Qin2022} associated with significant computational cost. In a number of tasks, the weighing is done not only for these different terms, but also for various collocation points chosen in the domain of interest~\cite{Wight2020}. Such weighting is based either on the minimization of the general loss function, or based on the features of some behavior, in particular, on the causality principle~\cite{Wang2022}, which has been efficiently used for a number of tasks and is actively leveraged to modify other methods and accelerate their convergence~\cite{Daw2022}. Notwithstanding the successes achieved, the PINN approach still continues to evolve and its development requires clear and transparent methods for weightings of coefficients of the elements of the general loss function and the architecture of the neural network for the problem.

This work is devoted to the development and application a number of techniques that can improve the accuracy of neural networks trained on the basis of the PINN approach. We will briefly list the contributions made in the paper:
\begin{enumerate}
	\item Reformulation of the original problem described by a differential equation accompanied with initial and boundary conditions to the form in which takes advantage of the differential equation only.
	\item Based on the reformulated problem, the loss function with the single significant term is proposed.
	\item Further, the modified loss function based on the causality principle is derived to account the possibility of using it for large number of temporal slices.
	\item Loss functions based on causality and spatial-temporal locality principles are proposed within the framework of the formalism of generalized functions (Heaviside step function and Dirac delta function).
	\item A number of experiments have been carried out for the above techniques for fully connected neural networks with homogeneous and heterogeneous activation functions.
\end{enumerate}

The paper is structured as follows. In section 2, the original PINN method and suggestions for modifications of the problems and loss function are presented. Section 3 includes the original description of causal training and modifications of the causal training and its improvement to spatial-temporal training within the framework of the formalism of generalized functions. In Section 4, numerical experiments to illustrate the accuracy and the efficiency of presented methods are given. Finally, in the Section 5 concluding remarks are given.

\section{Formulation of the problem}

\subsection{Physical-informed neural networks (PINNs)}
Consider nonlinear partial differential equations (PDEs), which depend on time $t$ and spatial coordinates ${\vec x}$, and in general take the following form
\begin{equation}\label{eq1}
	{\vec u}_{t} + \mathcal{N} [{\vec u}] = 0, \quad t \in [0, T], \quad {\vec x} \in \Omega,
\end{equation}
under the initial and boundary conditions
\begin{eqnarray}
	&& {\vec u}(0, {\vec x}) = {\vec g}({\vec x}), \quad {\vec x} \in \Omega, \label{eq2}\\
	&& \mathcal{B}[{\vec u}] = 0, \quad t \in [0, T], \quad {\vec x} \in \partial\Omega, \label{eq3}
\end{eqnarray}
where ${\vec u}(t,{\vec x})$ denotes the latent solution that is governed by the PDE system of equation (\ref{eq1}), ${\vec u}_t$ is temporal derivative, $\mathcal{N}$ is a nonlinear differential operator, $\mathcal{B}$ is a boundary operator corresponding to boundary conditions (Dirichlet, Neumann, Robin, periodic conditions), ${\vec {g}}({\vec x})$ is initial distribution of ${\vec u}$, $\Omega$ and $\partial\Omega$ are spatial domain and its boundary, respectively.

According to PINN approach~\cite{Raissi2019}, which stands on the shoulders of the universal approximation theorem~\cite{HORNIK1989359}, the aim of PINN is to approximate the unknown solution ${\vec u}(t, {\vec x})$ with a deep neural network ${\vec u}_{{\bm \theta}}(t,{\vec x})$, where ${\bm \theta}$ denote all trainable parameters (e.g., weight and biases) of the network. Finding optimal parameters is an optimization problem, which requires definition of a loss function such that its minimum gives the solution of the PDE. The physical-informed model is trained by minimizing the composite loss function which consists of the local residuals of the differential equation over the problem domain and its boundary as shown below:

\begin{equation}\label{eq4}
	\mathcal{L}({\bm \theta}) = \lambda_{ic} \mathcal{L}_{ic}({\bm \theta}) + \lambda_{bc} \mathcal{L}_{bc}({\bm \theta}) + \lambda_{r} \mathcal{L}_{r}({\bm \theta}),	  
\end{equation}
where
\begin{eqnarray}
	&& \mathcal{L}_{ic} \left( {\bm \theta} \right)  = \frac{1}{N_{ic}} \sum_{i=1}^{N_{ic}} \left| {\vec u}_{{\bm \theta}}\left(0, {\vec x}_{i}^{(ic)}\right) - {\vec g}\left( {\vec x}_{i}^{(ic)} \right)  \right|^{2},\label{eq5}\\
	&& \mathcal{L}_{bc}\left(\bm \theta \right) = \frac{1}{N_{bc}} \sum_{i=1}^{N_{bc}} \left| \mathcal{B}\left[{\vec u}_{{\bm \theta}} \right] \left(t_{i}^{(bc)}, {\vec x}_{i}^{(bc)}  \right)  \right|^{2}\label{eq6},\\
	&& \mathcal{L}_{r}\left(\bm \theta \right) = \frac{1}{N_r} \sum_{i=1}^{N_{r}} \left| \mathcal{R}\left[{\vec u}_{{\bm \theta}} \right] \left(t_{i}^{(r)}, {\vec x}_{i}^{(r)} \right)\right|^{2},\label{eq7}\\
	&& \mathcal{R}\left[{\vec u} \right]:= {\vec u}_{t} + \mathcal{N}\left[ {\vec u}\left(t, {\vec x} \right)  \right].\label{eq8}
\end{eqnarray}
Here $\left\{{\vec x}_{i}^{(ic)}\right\}^{N_{ic}}_{i=1}$,  $\left\{{t}_{i}^{(bc)}, {\vec x}_{i}^{(bc)}\right\}^{N_{bc}}_{i=1}$, and $\left\{{t}_{i}^{(r)}, {\vec x}_{i}^{(r)}\right\}^{N_{r}}_{i=1}$ are sets of points corresponding to initial condition domain, boundary condition domain and PDE domain, respectively. These points can be the vertices of a fixed mesh or can be randomly sampled at each iteration of a gradient descent algorithm. All required gradients w.r.t. input variables ($t$ and $\vec x$) and network parameters $\theta$ can be efficiently computed via automatic differentiation~\cite{Griewank2008} with algorithmic accuracy, which is defined by the accuracy of computation system. Note, the hyperparameters $\lambda_{ic}$, $\lambda_{bc}$ and $\lambda_{r}$ allow for separate tuning of learning rate for each of the loss terms in order to improve convergence of the model~\cite{Wang2021,WANG2022110768}.

\begin{figure}[t!]
	\centering
	\includegraphics[width=1\textwidth]{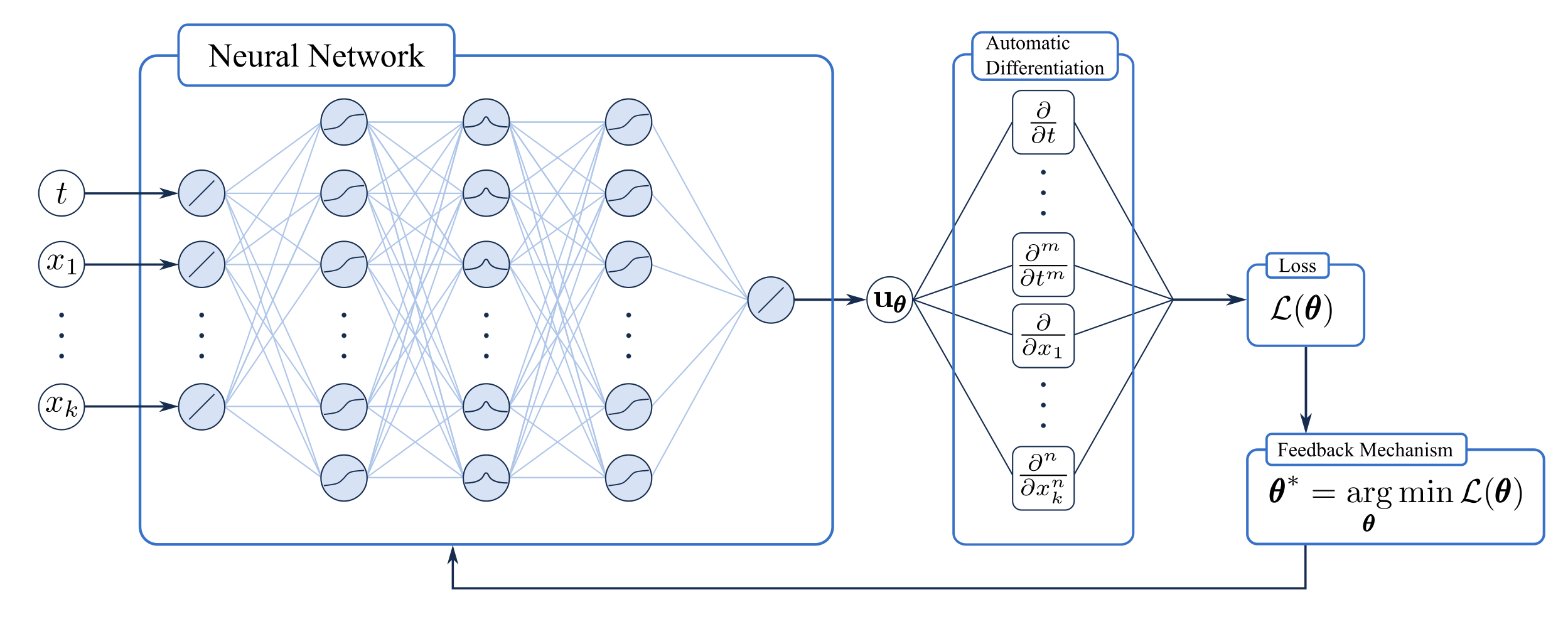}
	\caption{Schematic diagram of the general PINN method.}
	\label{fig1}
\end{figure}

The optimization problem can be defined as follow
\begin{equation}
	{\bm \theta}^* = \underset{{\bm \theta}}{\arg}\min \mathcal{L}({\bm \theta}),\label{eq9}
\end{equation}
where ${\bm \theta}^*$ are optimal parameters of the neural network which minimize the discrepancy between the exact unknown solution ${\vec u}$ and the approximate one ${\vec u}_{\pmb \theta^*}$. Schematic diagram of the general PINN method, which is used throughout this paper, is shown on the Figure~\ref{fig1}.

Let's rewrite the loss components (\ref{eq5})--(\ref{eq7}) in the continuous form for the further analysis. These representations are
\begin{align}
	& \mathcal{L}_{ic} \left( {\bm \theta} \right)  = \frac{1}{V_\Omega} \int_\Omega \left| {\vec u}_{{\bm \theta}}\left(0, {\vec x}\right) - {\vec g}\left( {\vec x} \right)  \right|^{2} d {\vec x},\label{eq10}\\
	& \mathcal{L}_{bc}\left(\bm \theta \right) = \frac{1}{T}\int_0^T \left| \mathcal{B}\left[{\vec u}_{{\bm \theta}} \right] \left(t, {\vec x}\right)  \right|^{2}d t,\quad {\vec x} \in {\partial \Omega}\label{eq11},\\
	& \mathcal{L}_{r}\left(\bm \theta \right) = \frac{1}{T}\frac{1}{V_\Omega}\int^{T}_0 dt \int_\Omega d{\vec x}  \left| \mathcal{R}\left[{\vec u}_{{\bm \theta}} \right] \left(t, {\vec x} \right)\right|^{2}\label{eq12},
\end{align}
where $V_\Omega$ is the ``volume'' of domain $\Omega$.

\subsection{Reformulation of the problem}

\subsubsection{Problem in time domain}
Let's modify the problem (\ref{eq1})--(\ref{eq3}) to the one which is free from the initial conditions. For this we expand the time interval from infinitely distant moment in the past to the time $T$. On time interval $(-\infty, 0]$ the value of given process ($\vec U$) is assumed to be constant in time and has a spatial distribution described by the initial distribution (\ref{eq2}). The value $\vec U$ has the following form:
\begin{equation}\label{eq13}
	{\vec U} = {\vec g}({\vec x})\chi(-t) + {\vec u}(t,{\vec x})\chi(t), \quad t\in (-\infty, T],
\end{equation}
where ${\vec u}(t,{\vec x})$ is unknown function that supposed to be a solution of the system of equations (\ref{eq1})--(\ref{eq3}), $\chi(t)$ is Heaviside step function~\cite{Davies2002} defined as
\begin{equation}\label{eq14}
	\chi(t) = \left\{\begin{array}{ll}
		0, & t < 0,\\
		1/2, & t=0,\\
		1, & t>0.
	\end{array}\right.
\end{equation}
In order the new function $\vec U$ to be the solution at time interval $(-\infty, 0]$ the initial equation (\ref{eq1}) should be supplied with additional term. Such term corresponds to ``force'' which disappears instantly at the moment $t=0$.

Without the loss of generality and for the sake of simplicity, we will consider the problem in an unbounded domain, i.e. without boundary conditions. As result the problem (\ref{eq1})--(\ref{eq3}) is reduced to the following differential equations
\begin{equation}\label{eq15}
	{\vec U}_{t} + \mathcal{N} [{\vec U}] = \mathcal{N} [{\vec g}({\vec x})]\chi(-t), \quad t \in (-\infty, T], \quad {\vec x} \in \Omega,
\end{equation}
where unknown function ${\vec u}(t,{\vec x})$ is encapsulated into $\vec U$. As previously, we approximate the unknown solution ${\vec U}(t, {\vec x})$ with a deep neural network ${\vec U}_{{\bm \theta}}(t,{\vec x}) = {\vec g}({\vec x})\chi(-t) + {\vec u}_{{\bm \theta}}(t,{\vec x})\chi(t)$.

In this problem statement, the loss function for problem (\ref{eq15}) contains only the term corresponding to residuals of differential equation:
\begin{align}
	& \mathcal{L}({\bm \theta}) = \mathcal{L}_{r}({\bm \theta}),\label{eq16}
\end{align}
where
\begin{align}
	& \mathcal{L}_{r}\left(\bm \theta \right) = \frac{1}{T}\frac{1}{V_\Omega}\int^{T}_{-\infty} dt \int_\Omega d{\vec x}  \left| \mathcal{R}^*\left[{\vec U}_{{\bm \theta}}, {\vec g}\right] \left(t, {\vec x} \right)\right|^{2},\label{eq17}\\
	& \mathcal{R}^*\left[{\vec u}, {\vec g}\right]:= {\vec u}_{t} + \mathcal{N}\left[ {\vec u}\left(t, {\vec x} \right)\right] - \mathcal{N} [{\vec g}({\vec x})]\chi(-t).\label{eq18}
\end{align}
Expanding the integrand in (\ref{eq17}) we have
\begin{equation}
	\hspace{-5mm}\mathcal{R}^*\left[{\vec U}_{{\bm \theta}}, {\vec g}\right] \left(t, {\vec x} \right) {=} \left\{\begin{array}{ll}
		\mathcal{R}\left[{\vec u}_{{\bm \theta}} \right] \left(t, {\vec x} \right) {+} \left[{\vec u}_{{\bm \theta}}\left(t, {\vec x} \right) {-} {\vec g}\left({\vec x} \right) \right] \delta(t), & t\geq 0,\\
		0 & t<0.
	\end{array}\right. \label{eq19}
\end{equation}
Here $\delta(t)$ is Dirac delta function~\cite{Saichev1996}. We use the following properties of delta function and Heaviside step function~\cite{Saichev1996}
\begin{eqnarray}
	&& \delta(x) = \frac{\partial \chi(x)}{\partial x}, \quad \chi(-x) = 1 - \chi(x). \label{eq20}
\end{eqnarray}
Substitution of (\ref{eq19}) to (\ref{eq17}) gives expression for loss in the following representation
\begin{align}
	\mathcal{L}_{r}\left(\bm \theta \right) & = \frac{1}{T}\frac{1}{V_\Omega}\lim_{t \to 0} \delta(t)\int_\Omega d{\vec x}  \left| {\vec u}_{{\bm \theta}}\left(0, {\vec x} \right) {-} {\vec g}\left({\vec x} \right)\right|^{2} \notag\\
	& + 2 \frac{1}{T}\frac{1}{V_\Omega} \int_\Omega d{\vec x}  \left| \mathcal{R}\left[{\vec u}_{{\bm \theta}} \right] \left(0, {\vec x} \right) \left[{\vec u}_{{\bm \theta}}\left(0, {\vec x} \right) {-} {\vec g}\left({\vec x} \right) \right]\right|\notag\\
	& + \frac{1}{T}\frac{1}{V_\Omega}\int^{T}_0 dt \int_\Omega d{\vec x}  \left| \mathcal{R}\left[{\vec u}_{{\bm \theta}} \right] \left(t, {\vec x} \right)\right|^{2}.\label{eq21}
\end{align}
Here we have taken into account the property of delta function which states its dominance in the finite sum, i.e. $\left|a \delta(x) + b\right|^2 \simeq \left|a \right|^2\delta^2(x) + 2 \left|a b\right| \delta(x) + \left|b\right|^2$.

\subsection{Optimization problem reformulation}

Then we apply the well known representation of the Dirac delta function~\cite{Saichev1996}
\begin{align}
	& \delta(x) = \lim_{\sigma \to 0} \frac{\sigma}{x^2 + \sigma^2}.\label{eq22} 
\end{align}
We obtain the following formulation of the loss function for (\ref{eq9})
\begin{align}
	\mathcal{L}\left(\bm \theta \right) & = \frac{1}{T}\frac{1}{V_\Omega}\int^{T}_0 dt \int_\Omega d{\vec x}  \left| \mathcal{R}\left[{\vec u}_{{\bm \theta}} \right] \left(t, {\vec x} \right)\right|^{2} \notag\\
	& +\frac{1}{T}\frac{1}{V_\Omega}\left\{ 2 \int_\Omega d{\vec x}  \left| \mathcal{R}\left[{\vec u}_{{\bm \theta}} \right] \left(0, {\vec x} \right) \left[{\vec u}_{{\bm \theta}}\left(0, {\vec x} \right) {-} {\vec g}\left({\vec x} \right) \right]\right|\right.\notag\\
	&\left. +\lim_{\substack{ \sigma \to 0 \\ t \to 0}} \frac{\sigma}{t^2 + \sigma^2}\int_\Omega d{\vec x}  \left| {\vec u}_{{\bm \theta}}\left(0, {\vec x} \right) {-} {\vec g}\left({\vec x} \right)\right|^{2} \right\}.\label{eq23}
\end{align}

To proceed from continuous time to the discrete time representation, we divide the time interval $[0,T]$ into $N_t$ intervals and use the rectangular rule of integration for this mesh with step size $T/N_t$. We take $\sigma = T/N_t$ for the approximation of the delta function in the time domain. This representation of delta function is valid as it preserves its fundamental property~\cite{Saichev1996} [$\int_{-\infty}^{\infty}f(x)\delta(x-a)dx = f(a)$] for the numerical integration. Finally, we get the following loss function	 
\begin{align}
	\mathcal{L}\left(\bm \theta \right) & = \frac{1}{N_t}\frac{1}{V_\Omega}\sum\limits_{i=0}^{N_t} \int_\Omega d{\vec x}  \left| \mathcal{R}\left[{\vec u}_{{\bm \theta}} \right] \left(t_i, {\vec x} \right)\right|^{2} \notag\\
	& + \frac{1}{T}\frac{1}{V_\Omega}\left\{2 \int_\Omega d{\vec x}  \left| \mathcal{R}\left[{\vec u}_{{\bm \theta}} \right] \left(t_0, {\vec x} \right) \left[{\vec u}_{{\bm \theta}}\left(t_0, {\vec x} \right) {-} {\vec g}\left({\vec x} \right) \right]\right|\right.\notag\\
	&\left. + \frac{N_t}{T} \int_\Omega d{\vec x}  \left| {\vec u}_{{\bm \theta}}\left(t_0, {\vec x} \right) {-} {\vec g}\left({\vec x} \right)\right|^{2} \right\}.\label{eq24}
\end{align}
Here we take $t_0 = 0$. Note, the multiplier ${N_t}$ of third term of (\ref{eq24}) coincides with coefficient $\lambda_{ic}$ empirically found in ~\cite{Wang2022,Wang2022_2,Zhang2023} in the case of $N_t=100$.

With this representation, the approximation of delta function is associated with a numerical integration scheme. For a number of tasks though, an approximation of the delta function is possible that is not tied to integration rules. To achieve that the value $\sigma$ in (\ref{eq22}) can be taken as an estimated relaxation time of the initial conditions to zero, or of doubling or halving values of initial conditions. Note that obtained approximation of the delta function can be used as a weight $\lambda_{ic}$ in the original training scheme described by loss terms~(\ref{eq5})--(\ref{eq7}). Example of deriving such approximation for specific problem is provided in section~\ref{approxDelta}.

\subsubsection{MAE loss}
Interestingly, if instead of MSE loss we use MAE loss for optimization problem (\ref{eq9}), then no additional representation of the Dirac delta function is required and instead of loss~(\ref{eq23}) we have the following function
\begin{align}
	\mathcal{L}\left(\bm \theta \right) & = \frac{1}{T}\frac{1}{V_\Omega} \left( \int^{T}_0 dt \int_\Omega d{\vec x}  \left| \mathcal{R}\left[{\vec u}_{{\bm \theta}} \right] \left(t, {\vec x} \right)\right|  + \int_\Omega d{\vec x}  \left|{\vec u}_{{\bm \theta}}\left(0, {\vec x} \right) {-} {\vec g}\left({\vec x} \right)\right|\right).\label{eq25}
\end{align}
Based on our experiments, training a neural network with such a loss function converges extremely weak to the required solution, despite the simplicity of the equation~(\ref{eq25}). For this reason, in all of the experiments presented in section~\ref{experiments}, MSE loss has been used.

\subsection{Problem in spatial-temporal domain}\label{appB}
We can reformulate the problem in spatial domain in the same manner. For greater clarity and brevity, and without violating generality, we will consider a problem with dependence on only one spatial coordinate $x$. 

Let spatial domain be determined by a coordinate $x$ ($x \in [X_1, X_2] $), and problem (\ref{eq1})--(\ref{eq3}) is subject to the following boundary conditions
\begin{align}
	& {\vec u}(t, X_1) = {\vec f_1}(t), \quad {\vec u}(t, X_2) = {\vec f_2}(t). \label{eq26}
\end{align}
Then, artificially extended solution of (\ref{eq1})--(\ref{eq3}) takes the following form
\begin{align}\label{eq27}
	{\vec U} &= {\vec g}({x})\chi(-t) + {\vec u}(t,{x})\chi(t)\left[\chi(x - X_1) - \chi(x - X_2)\right]\notag\\
	& + \left[{\vec f_1}(t)\chi(X_1 - x) + {\vec f_2}(t)\chi(x - X_2)\right]\chi(t),\\
	& {\text{for the domain }} t\in (-\infty, T],\quad x \in (-\infty, \infty). \notag 
\end{align}

In this case, the problem (\ref{eq1})--(\ref{eq3}) is reduced to the following differential equations
\begin{align}\label{eq28}
	{\vec U}_{t} + \mathcal{N} [{\vec U}] &= \mathcal{N} [{\vec g}({x})]\chi(-t) \notag\\ &+ \left[\frac{\partial {\vec f_1}(t)}{\partial t}\chi(X_1 - x)
	+ \frac{\partial{\vec f_2}(t)}{\partial t} \chi(x - X_2)\right]\chi(t)\notag\\
	&+ \left\{ \mathcal{N}\left[{\vec f_1}(t)\right]\chi(x - X_1)
	+ \mathcal{N}\left[{\vec f_2}(t)\right] \chi(X_2 - x)\right\}\chi(t), \quad t\in (-\infty, T],\quad x \in (-\infty, \infty).
\end{align}
The value $\mathcal{R}^*\left[{\vec U}_{{\bm \theta}}, {\vec g}\right] \left(t, {x} \right)$ has form similar to (\ref{eq19})
\begin{equation}
	\hspace{-5mm}\mathcal{R}^*\left[{\vec U}_{{\bm \theta}}, {\vec g}\right] \left(t, {x} \right) {=} \left\{\begin{array}{lll}
			\mathcal{R}\left[{\vec U}^*_{{\bm \theta}} \right] \left(t, {x} \right)&\hspace{-3mm}{+} \left[{\vec u}_{{\bm \theta}}\left(t, {x} \right) {-} {\vec g}\left({x} \right) \right] \delta(t) & \\
			&\hspace{-13mm}- \left[\frac{\partial {\vec f_1}(t)}{\partial t}\chi(X_1 - x)
			+ \frac{\partial{\vec f_2}(t)}{\partial t} \chi(x - X_2)\right] & \\
			&\hspace{-13mm}- \left\{ \mathcal{N}\left[{\vec f_1}(t)\right]\chi(x - X_1)
			+ \mathcal{N}\left[{\vec f_2}(t)\right] \chi(X_2 - x)\right\}
		 , & t\geq 0, \quad x \in [X_1, X_2],\\
		0, & & t<0,\quad \text{or } x \notin [X_1, X_2],
	\end{array}\right. \label{eq29}
\end{equation}
where 
\begin{align}\label{eq30}
	{\vec U}^*_{{\bm \theta}} &= {\vec u}_{{\bm \theta}}(t,{x})\left[\chi(x - X_1) - \chi(x - X_2)\right] + {\vec f_1}(t)\chi(X_1 - x) + {\vec f_2}(t)\chi(x - X_2).
\end{align}
The representation of $\mathcal{R}\left[{\vec U}^*_{{\bm \theta}} \right]$ in terms of $\mathcal{R}\left[{\vec u}_{{\bm \theta}} \right]$ similarly to equation~(\ref{eq19}) depends on the specific type of differential equation. Example of such representation can be found in section~\ref{KdV}.

\section{Causal training for physics-informed neural network}

It was suggested in~\cite{Wang2022} to reformulate the loss function so that it would respect the physical causality when solving PDEs with the help of PINNs:
\begin{equation}
	\mathcal{L}_{r}\left(\bm \theta \right) = \frac{1}{N_t}\sum_{i=0}^{N_{t}} {w}_{i} \mathcal{L}_{r}\left(t_{i}, {\bm \theta} \right),\label{eq31}  
\end{equation}
where weights $w_i$ are
\begin{equation}
	w_{i} = \exp \left[ - \varepsilon \sum_{k=0}^{i-1} \mathcal{L}_r \left( t_k, {\bm \theta} \right)  \right], \quad \text{for } i=1,2,3,\ldots,N_{t}, \text{and } w_0 =1,\label{eq32}
\end{equation}
and where $\varepsilon$ is causality parameter that controls the steepness of the weights $w_i$. The temporal residual loss [$\mathcal{L}_{r}\left(t_{i}, {\bm \theta} \right)$] for the total loss (\ref{eq23}) is
\begin{align}
	\mathcal{L}_{r}\left(t_i, {\bm \theta} \right) = &\frac{1}{V_\Omega} \int_\Omega d{\vec x} \left\{ \left| \mathcal{R}\left[{\vec u}_{{\bm \theta}} \right] \left(t_i, {\vec x} \right)\right|^{2} \right.\notag\\
	& + 2 \frac{N_t}{T}\delta_{0i} \left| \mathcal{R}\left[{\vec u}_{{\bm \theta}} \right] \left(t_0, {\vec x} \right) \left[{\vec u}_{{\bm \theta}}\left(t_0, {\vec x} \right) {-} {\vec g}\left({\vec x} \right) \right]\right|\notag\\
	&\left. + \frac{N_t^2}{T^2}\delta_{0i} \left| {\vec u}_{{\bm \theta}}\left(t_0, {\vec x} \right) {-} {\vec g}\left({\vec x} \right)\right|^{2} \right\}, \label{eq33}
\end{align}
$\delta_{0i}$ is the Kronecker symbol.

In the case when domain $\Omega$ consists of single spatial dimension $x$ ($x \in [X_1, X_2]$, $V_\Omega = X_2 - X_1$) the temporal residual loss has discrete form under regular mesh with $N_x$th nodes
\begin{align}
	\mathcal{L}_{r}\left(t_i, {\bm \theta} \right) = & \frac{1}{N_x} \sum_{n=0}^{N_x} \left\{ \left| \mathcal{R}\left[{\vec u}_{{\bm \theta}} \right] \left(t_i, {x_n} \right)\right|^{2} \right.\notag\\
	& + 2 \frac{N_t}{T}\delta_{0i} \left| \mathcal{R}\left[{\vec u}_{{\bm \theta}} \right] \left(t_0, {x_n} \right) \left[{\vec u}_{{\bm \theta}}\left(t_0, {x_n} \right) {-} {\vec g}\left({x_n} \right) \right]\right|\notag\\
	&\left. + \frac{N_t^2}{T^2}\delta_{0i} \left| {\vec u}_{{\bm \theta}}\left(t_0, {x_n} \right) {-} {\vec g}\left({x_n} \right)\right|^{2} \right\}. \label{eq34}
\end{align}
The resulting weighted residual loss is
\begin{equation}
	\mathcal{L}_{r}\left(\bm \theta \right) = \frac{1}{N_t}\sum_{i=0}^{N_{t}} \exp \left[ - \varepsilon \sum_{k=0}^{i-1} \mathcal{L}_r \left( t_k, {\bm \theta} \right)  \right] \mathcal{L}_{r} \left(t_{i}, {\bm \theta} \right).\label{eq35}  
\end{equation}
Authors of~\cite{Wang2022} recognized that the weights $w_i$ should be large --- and therefore allow the minimization of $\mathcal{L}_{r}\left( t_i, {\bm \theta} \right) $ --- only if all residuals $\left\{ \mathcal{L}_{r}\left( t_k, {\bm \theta} \right)  \right\}_{k=0}^{i}$ before $t_i$ are minimized properly, and vice versa. As a consequence, $\mathcal{L}_{r}\left(t_{i}, {\bm \theta} \right) $ is not minimized unless all previous residuals $\left\{ \mathcal{L}_{r}\left( t_k, {\bm \theta} \right)  \right\}_{k=0}^{i-1}$ decrease to some small value such that $w_i$ is large enough.

\subsection{Modification of the causal training for physics-informed neural network}

Despite the causal training method significantly improves the convergence of neural network learning this approach does not scale for time grids with the growth of number of nodes. This can be explained by the fact that the exponent of the weighting function (\ref{eq32}) contains the cumulative sum. This sum is increasing with the increase of the $N_t$ and such behavior hinders learning of neural network for large $t$. Note that this feature in row of cases can be lead to the best convergence of training in compared to other methods of weighting.

This drawback becomes obvious when trying to write the residual loss (\ref{eq35}) in continuous form
\begin{equation}
	\mathcal{L}_{r}\left(\bm \theta \right) = \frac{1}{T} \int_{0}^{T}  \lim_{N_t\to \infty}\left[\exp \left( - \varepsilon\frac{N_t}{t} \int_{0}^{t} \mathcal{L}_r \left( t', {\bm \theta} \right) dt'  \right)\right] \mathcal{L}_{r} \left(t, {\bm \theta} \right)d t.	\label{eq36}
\end{equation}
This integral tends to zero for any $\varepsilon$ and time due to unlimited increase of $N_t$. Such disadvantage prevents mathematical analysis of the loss function using standard methods based on integral and differential calculus.

In order to mitigate the above drawback, the residual loss (\ref{eq35}) can be rewritten in the following form
\begin{align}
	\mathcal{L}_{r}\left(\bm \theta \right) & = \frac{1}{N_t}\sum_{i=0}^{N_{t}} w_i \mathcal{L}_{r} \left(t_{i}, {\bm \theta} \right),\label{eq37} \\
	w_i & = \exp \left[ - \varepsilon s_i  \right],\label{eq38} \\
	s_i & = \frac{1}{i} \sum_{k=0}^{i-1} \mathcal{L}_r \left( t_k, {\bm \theta} \right). \label{eq39}
\end{align}
The continuous form of residual loss (\ref{eq37}) is non zero and has the following representation
\begin{equation}
	\mathcal{L}_{r}\left(\bm \theta \right) = \frac{1}{T} \int_{0}^{T}d t  \exp \left[ - \frac{\varepsilon}{t} \int_{0}^{t} dt' \mathcal{L}_r \left( t', {\bm \theta} \right)  \right] \mathcal{L}_{r} \left(t, {\bm \theta} \right).	\label{eq40}
\end{equation}

\subsection{Alternatives of the causal training for physics-informed neural network}

\subsection{Heaviside--weighted residual loss}
Let's consider residual loss weighted using Heaviside step function
\begin{equation}
	\mathcal{L}_{r}\left(\bm \theta \right) = \frac{1}{T}\int_{0}^{T} \left\{1-
	\chi \left[S(t)\right]\right\}  \mathcal{L}_{r} \left(t, {\bm \theta} \right)d t, \label{eq41}
\end{equation}
where
\begin{equation}
	S(t) = \frac{1}{t} \int_{0}^{t} \mathcal{L}_r \left(t', {\bm \theta} \right)dt'. \label{eq42}
\end{equation}
In this case term $1-\chi(x)$ in the integrand of equation (\ref{eq41}) does not tend to zero under condition $x\leq 0$. It means that causality condition is satisfied. To calculate the loss we use the following well known representations of the Heaviside step function~\cite{Saichev1996}
\begin{equation}
	\chi(x) = \lim_{\varepsilon \to \infty} \frac{1}{1+ \exp\left(-2 \varepsilon x\right)}.\label{eq43}
\end{equation}
The equation (\ref{eq43}) leads to the following representation of residual loss
\begin{equation}
	\mathcal{L}_{r}\left(\bm \theta \right) = \frac{1}{T}\int_{0}^{T}d t  \frac{1}{1+ \exp\left[2 \varepsilon S(t)\right]} \mathcal{L}_{r} \left(t, {\bm \theta} \right),	\label{eq44}
\end{equation}
which has the discrete form (\ref{eq31}) with weights
\begin{align}
	& w_i = \left[1 + \exp  \left({ 2 \varepsilon} s_i\right) \right]^{-1} \label{eq45}.
\end{align}

\subsubsection{Dirac--weighted residual loss}

Let's consider the extreme form of the weighted residual loss
\begin{equation}
	\mathcal{L}_{r}\left(\bm \theta \right) = \frac{1}{T}\int_{0}^{T} \delta \left[S(t)\right]  \mathcal{L}_{r} \left(t, {\bm \theta} \right)d t. \label{eq46}
\end{equation}
Such loss function has well-pronounced property that it prohibits training for large time stamps without have been previously trained for small ones.

To proceed with the discrete form of this loss function we will use the following well known representation of the Dirac delta function~\cite{Saichev1996}
\begin{align}
	& \delta(x) = \lim_{\sigma^2 \to 0} \frac{1}{\sigma\sqrt{2\pi}}\exp\left( - \frac{x^2}{2 \sigma^2} \right).\label{eq47}
\end{align}

The integral form of residual loss in this case is
\begin{equation}
	\mathcal{L}_{r}\left(\bm \theta \right) = \frac{1}{T}\int_{0}^{T}d t  \exp \left\{ - \varepsilon \left[\frac{1}{t} \int_{0}^{t} dt' t' \mathcal{L}_r \left( t', {\bm \theta} \right)\right]^2  \right\} \mathcal{L}_{r} \left(t, {\bm \theta} \right).	\label{eq48}
\end{equation}
Note, here the causality properties of weights are further enhanced with multiplying by $t$ inside the inner integral. The discrete form of (\ref{eq48}) in case of equidistant time slices is
\begin{align}
	& \mathcal{L}_{r}\left(\bm \theta \right) = \frac{1}{N_t}\sum_{i=0}^{N_{t}} w_i \mathcal{L}_{r} \left(t_{i}, {\bm \theta} \right),\label{eq49} \\
	& w_i = \exp \left(- \varepsilon \tilde{s}_i^2  \right), \label{eq50}\\
	& \tilde{s}_i = \frac{1}{i} \sum_{k=0}^{i-1} k \mathcal{L}_r \left( t_k, {\bm \theta} \right). \label{eq51}
\end{align}
When deriving the equation~(\ref{eq51}) we have taken into account that $t_k = k T / N_t$ and multiplier $T / N_t$ of $\tilde{s}_i$ was ``hidden'' in causality parameter $\varepsilon$ in (\ref{eq50}).

\subsection{Expanding of the causal training approach to the spatial-temporal domain}
For problem (\ref{eq28}) the loss function has the following representation
\begin{align}
	\mathcal{L}\left(\bm \theta \right) & = \frac{1}{T}\frac{1}{V_\Omega}\int^{T}_0 dt \int_\Omega d{x}  \left| \mathcal{R}^*\left[{\vec U}^*_{{\bm \theta}} \right] \left(t, {x} \right)\right|^{2} \notag\\
	& +\frac{1}{T}\frac{1}{V_\Omega}\left\{ 2 \int_\Omega d{x}  \left| \mathcal{R}^*\left[{\vec U}^*_{{\bm \theta}} \right] \left(0, {x} \right) \left[{\vec u}_{{\bm \theta}}\left(0, {x} \right) {-} {\vec g}\left({x} \right) \right]\right|\right.\notag\\
	&\left. +\lim_{\substack{ \sigma \to 0 \\ t \to 0}} \frac{\sigma}{t^2 + \sigma^2}\int_\Omega d{x}  \left| {\vec u}_{{\bm \theta}}\left(0, {x} \right) {-} {\vec g}\left({x} \right)\right|^{2} \right\}.\label{eq52}
\end{align}
The weighted residual loss in the spatial-temporal domain can be defined as follows
\begin{align}
	& \mathcal{L}^*_{r}\left(\bm \theta \right) = \frac{1}{N_t}\sum_{i=0}^{N_{t}} w_i \mathcal{L}^*_{r} \left(t_{i}, {\bm \theta} \right),\label{eq53} \\
	& w_i = \exp \left[ - \varepsilon s_i  \right],\quad s_i = \frac{1}{i} \sum_{k=0}^{i-1} \mathcal{L}^*_r \left( t_k, {\bm \theta} \right), \label{eq54}
\end{align}
where 
\begin{align}
	\mathcal{L}^*_{r}\left(t_i, {\bm \theta} \right) & = \frac{1}{N_x} \sum_{n=0}^{N_x} \left(w_n^{(\to)} + w_n^{(\leftarrow)}\right)\left\{ \left| \mathcal{R}\left[{\vec U}^*_{{\bm \theta}} \right] \left(t_i, {x_n} \right)\right|^{2} \right.\notag\\
	& + 2 \frac{N_t}{T}\delta_{0i} \left| \mathcal{R}\left[{\vec U}^*_{{\bm \theta}} \right] \left(t_0, {x_n} \right) \left[{\vec u}_{{\bm \theta}}\left(t_0, {x_n} \right) {-} {\vec g}\left({x_n} \right) \right]\right|\notag\\
	&\left. + \frac{N_t^2}{T^2}\delta_{0i} \left| {\vec u}_{{\bm \theta}}\left(t_0, {x_n} \right) {-} {\vec g}\left({x_n} \right)\right|^{2} \right\},\label{eq55}
\end{align}
and
\begin{align}
	w_n^{(\to)} &= \exp \left[ - \varepsilon s^{(\to)}_n  \right],\quad w_n^{(\leftarrow)} = \exp \left[ - \varepsilon s^{(\leftarrow)}_n  \right],\notag\\
	s^{(\to)}_n & = \frac{1}{n} \sum_{k=0}^{n-1} \tilde{\mathcal{L}}_r \left( x_k, {\bm \theta} \right),\quad  s^{(\leftarrow)}_n = \frac{1}{N_x - n} \sum_{k=0}^{N_x -n-1} \tilde{\mathcal{L}}_r \left( x_{N_x-k}, {\bm \theta} \right).\label{eq56}
\end{align}
The spatial residual loss is
\begin{align}
	\tilde{\mathcal{L}}_r \left( x_k, {\bm \theta} \right) & = \frac{1}{N_t} \sum_{i=0}^{N_t} \left| \mathcal{R}\left[{\vec u}_{{\bm \theta}} \right] \left(t_i, {x_k} \right)\right|^{2} . \label{eq57}
\end{align}
Here superscripts $(\to)$ and $(\leftarrow)$ mean weighting in the forward and backward directions along the $x$ axis, respectively (i.e. from $X_1$ to $X_2$ and from $X_2$ to $X_1$).

\subsection{Training procedure}

Based on the above modifications of the problem, loss function and causality weights, algorithms~\ref{alg1} and ~\ref{alg2} present a causal and a spatial-temporal local training algorithms for PINNs, respectively. 
\begin{algorithm}
	\caption{Training of physics-informed neural networks at time domain}\label{alg1}
	\KwData{Sets of coordinates $(t,{\vec x})$ and values ${\vec u}$ of points corresponding to initial and boundary condition domains and set of coordinates of collocation points}
	\KwResult{Trained NN ${\vec u}_{\pmb \theta}$ which approximates the solution of PDE}
	Initialize temporal weights $w_i$ with 1, causal parameter $\varepsilon$ and threshold parameter $\delta_w$ with chosen values.\\
	Use $N$ epochs of a gradient descent algorithm to update the parameters $\pmb \theta$ as follows \\
	\For{$n=1,\dots, N$}{
		1. Calculate and update the temporal weight $w_i$ given by equations (\ref{eq38}), (\ref{eq45}) or (\ref{eq50}) using the temporal residual loss (\ref{eq34})\\
		\If{$\underset{i}{\min} (\exp \left[ - \varepsilon s_i  \right]) > \delta_w$}{
			\begin{equation}
				\varepsilon \gets m_\varepsilon \times \varepsilon \notag
			\end{equation}
		}
		2. Calculate the loss using equations (\ref{eq16}) and (\ref{eq31})\\
		3. Update the parameters $\theta$ via gradient descent\\
		\begin{equation}
			{\bm \theta}_{n} \gets {\bm \theta}_{n-1} - \eta \nabla_{{\bm \theta}} \mathcal{L}\left(\bm \theta_{n-1} \right). \label{eq58}
		\end{equation}
	}
	From numerical experiments we found that the most promising hyperparameters are $\delta_w = 0.99$, initial value of $\varepsilon$=$10^{-3}$, multiplication factor $m_\varepsilon = 2$. Here $\eta$ is the learning rate of gradient descent.
\end{algorithm}

\begin{algorithm}
	\caption{Training of the physics-informed neural networks at spatio-temporal domain}\label{alg2}
	\KwData{Sets of coordinates $(t,{\vec x})$ and values ${\vec u}$ of points corresponding to initial and boundary condition domains and set of coordinates of collocation points}
	\KwResult{Trained NN ${\vec u}_{\pmb \theta}$ which approximates the solution of PDE}
	Initialize temporal and spatial weights $w_i$, $w_i^{(\to)}$, and $w_i^{(\leftarrow)}$ with 1, causal parameter $\varepsilon$ and threshold parameter $\delta_w$ with chosen values.\\
	Use $N$ epochs of a gradient descent algorithm to update the parameters $\pmb \theta$ as follows \\
	\For{$n=1,\dots, N$}{
		1. Calculate and update the spatial weights $w_i^{(\to)}$ and $w_i^{(\leftarrow)}$ given by equation (\ref{eq56}) using the spatial residual loss (\ref{eq57})\\
		2. Calculate and update the temporal weight $w_i$ given by equations (\ref{eq38}), (\ref{eq45}) or (\ref{eq50}) using the temporal residual loss (\ref{eq55})\\
		\If{$\underset{i}{\min} (\exp \left[ - \varepsilon s_i  \right]) > \delta_w$ \& $\underset{i}{\min} \left(\exp \left[ - \varepsilon s_i^{(\to)}  \right]\right) > \delta_w$ \& $\underset{i}{\min} \left(\exp \left[ - \varepsilon s_i^{(\leftarrow)}  \right]\right) > \delta_w$}{
			\begin{equation}
				\varepsilon \gets m_\varepsilon \times \varepsilon \notag
			\end{equation}
		}
		3. Calculate the loss using equations (\ref{eq16}) and (\ref{eq54})\\
		4. Update the parameters $\theta$ via gradient descent\\
		\begin{equation}
			{\bm \theta}_{n} \gets {\bm \theta}_{n-1} - \eta \nabla_{{\bm \theta}} \mathcal{L}\left(\bm \theta_{n-1} \right). \label{eq59}
		\end{equation}
	}
\end{algorithm}

Some details of the training procedures are explained below.

To save computational resources and to accelerate convergence for problems with periodic boundary conditions along $x$ axis let's construct a Fourier feature embedding of input coordinates in the following form
\begin{equation}\label{eq60}
	v(x) = \left[1,\cos(k_x x),\sin(k_x x),\cos(2 k_x x),\sin(2 k_x x),...,\cos(m k_x x),\sin(m k_x x)\right],
\end{equation}
where $k_x = 2\pi/L$, $L$ is the minimal spatial period along $x$ direction, and $m$ is some non-negative integer. With such embedding any trained approximation of problem solution ${\bm u}_{\pmb \theta}\left[v(x)\right]$ exactly satisfies periodic conditions (see details in~\cite{Wang2022,DONG2021110242}).

The networks had been trained via stochastic gradient descent using the
Adam optimizer~\cite{Kingma2014} and different schedulers (see details in the table~5).

\section{Numerical experiments}\label{experiments}

In this section several numerical experiments are presented to demonstrate the accuracy and performance aforementioned methods.

To evaluate the accuracy of the approximate solution obtained with the help of PINN method, the values of the solution of (\ref{eq1}) predicted by the neural network at given points are compared with the values calculated on the basis of classical high-precision numerical methods. As a measure of accuracy, the relative total ${\mathbb L}_2$ error of prediction is taken, which can be expressed with the following relation
\begin{equation}\label{eq61}
	\epsilon_{\rm error} = \left\{\frac{1}{N_e} \sum^{N_e}_{i=1} \left[{\vec u}_{\pmb \theta}(t_i, {\vec x}_i) - {\vec u}(t_i, {\vec x}_i)\right]^2 \right\}^{1/2} \times \left\{\frac{1}{N_e} \sum^{N_e}_{i=1} \left[{\vec u}(t_i, {\vec x}_i)\right]^2 \right\}^{-1/2},
\end{equation}
where $\left\{{t}_{i}, {\vec x}_{i}\right\}^{N_{e}}_{i=1}$ is the set of evaluation points taken from the domain $[0,T]\times \Omega$, ${\vec u}_{\pmb \theta}$ and ${\vec u}$ are the predicted and reference solutions respectively.

For our experiments we used Pytorch~\cite{Paszke2019} version 1.12.1 and the training was carried out on a node with Nvidia Tesla V100 GPU. 

\subsection{Allen--Cahn equation}
As in ~\cite{Raissi2019,Wang2022} let's consider a classical phase field modeled by one-dimensional Allen--Cahn equation with periodic boundary condition
\begin{align}
	& u_t - 0.0001 u_{xx} + 5 u^3 - 5u = 0, \quad t \in [0,1],\quad x\in[-1,1], \notag\\
	& u(x,0) = x^2 \cos(\pi x), \notag\\
	& u(t,-1) = u(t,1), \quad u_x(t,-1) = u_x(t,1).\label{eq62}
\end{align}
We represent the latent variable $u$ by a fully-connected neural network $u_{\pmb \theta}$ with hyperbolic tangent activation function, 4 hidden layers with 128 neurons per each layer. For simplicity, a uniform mesh of size $100\times256$ is constructed in computational domain $[0, 1]\times[-1, 1]$, yielding $N_{ic} = 256$ initial points and $N_r = 25600$ collocation points for enforcing the PDE residual. Fourier expansion of the input coordinates in the embedding is determined by the first 20 spatial harmonics ($m=20$). For this problem, the logarithm of (\ref{eq16}) was used as the loss function
for training $u_{\pmb \theta}$. We have chosen such an approach because based on our experiments it demonstrates faster and more stable convergence. The optimization problem in this case of loss function is
\begin{align}\label{eq63}
	{\bm \theta}^* = \underset{{\bm \theta}}{\arg}\min\log\left[{\mathcal{L}}(\pmb \theta)\right].
\end{align}

The reference solution can be generated using the Chebfun package~\cite{driscollChebfunGuide2014}, or dataset from supplementary materials for ~\cite{Wang2022} can be used.

\begin{figure}[t!]
	\centering
	\includegraphics[width=0.95\textwidth]{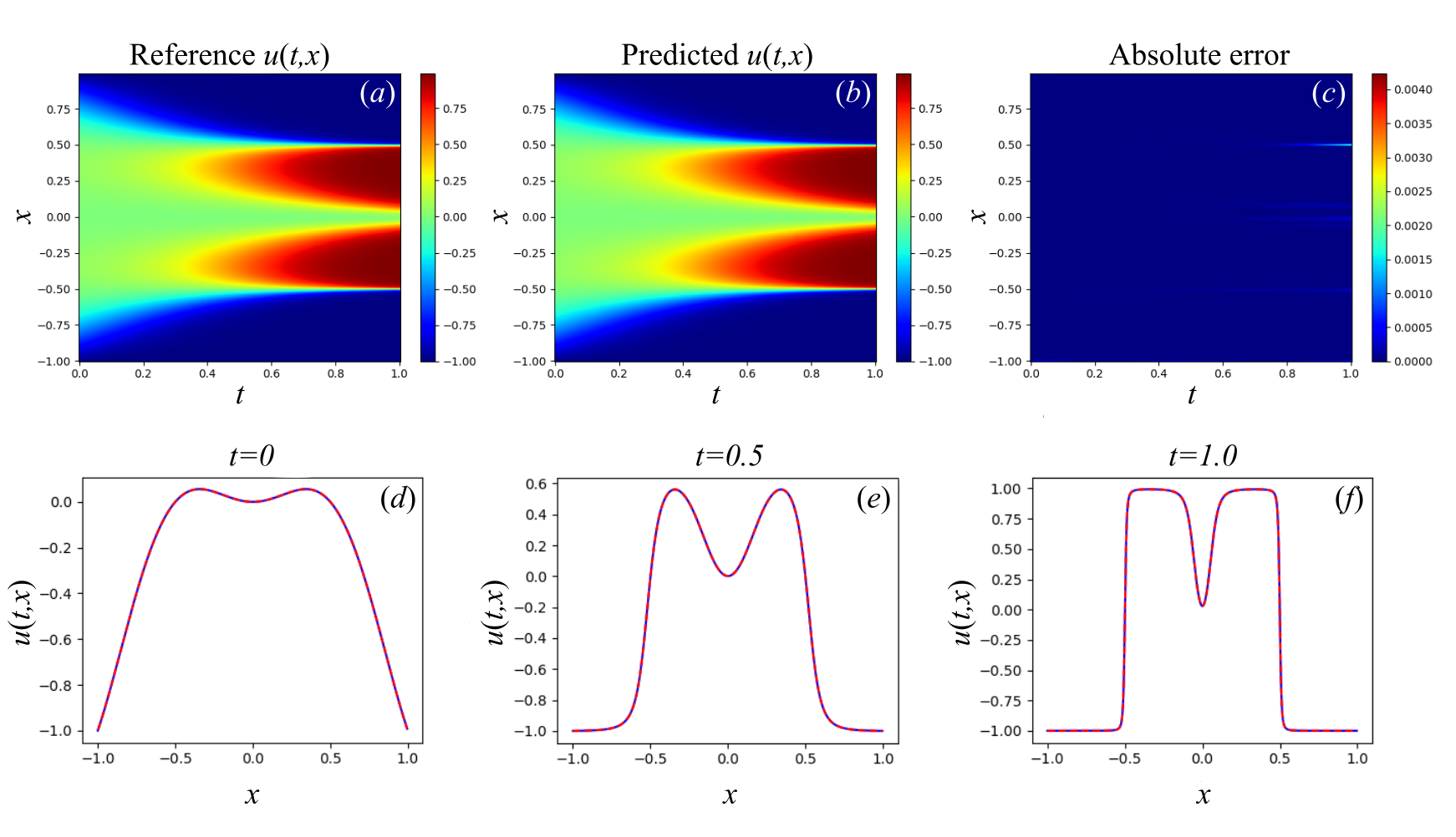}
	\caption{Allen--Cahn equation. (a) is reference solution, (b) is prediction of a trained physics-informed neural network based on algorithm~\ref{alg1}, (c) is absolute difference between reference solution and predicted solution. The relative error $\epsilon_{\rm error}$ is $6.29\times 10^{-5}$. (d), (e) and (f) are comparison of the predicted (red dash lines) and reference solutions (blue solid lines) corresponding to three temporal snapshots at $t = 0.0$, $t=0.5$ and $t=1.0$, respectively.}
	\label{fig2}
\end{figure}

The experiments were carried out for various weights described by equations~(\ref{eq38}), (\ref{eq45}) and (\ref{eq50}). Hyperparameters used for training are presented in table~5 and table~6.
In table~1 we also report the accuracy for Allen--Cahn problem achieved using other methods from the literature~\cite{Raissi2019,Wight2020,McClenny2020,Mattey2022,Wang2022,Wang2022_2,Zhang2023}. As can be seen from this table the best result is achieved with Dirac delta function and MLP with a relative ${\mathbb L}_2$ error $6.29\times 10^{-5}$. Results of this experiment are shown on the Figure~\ref{fig2}. The predicted solution coincides with the ground truth almost perfectly.

\begin{table}[!h]
	\begin{center}\label{table1}
		\begin{tabular}{l|c}
			\hline
			Method \rule[-1ex]{0pt}{3.5ex}  & Relative ${\mathbb L}_2$ error  \\
			\hline
			Original formulation by Raissi et al.~\cite{Raissi2019} \rule[-1ex]{0pt}{3.5ex}  & $4.98\times 10^{-1}$ \\
			Adaptive time sampling~\cite{Wight2020} \rule[-1ex]{0pt}{3.5ex}  & $2.33\times 10^{-2}$ \\
			Self-attention~\cite{McClenny2020} \rule[-1ex]{0pt}{3.5ex}  & $2.10\times 10^{-2}$ \\
			Time marching~\cite{Mattey2022} \rule[-1ex]{0pt}{3.5ex}  & $1.68\times 10^{-2}$ \\
			Causal training (MLP)~\cite{Wang2022} \rule[-1ex]{0pt}{3.5ex}  & $1.43\times 10^{-3}$ \\
			Modified MLP with conservation laws~\cite{Wang2022_2} \rule[-1ex]{0pt}{3.5ex}  & $5.42\times 10^{-4}$ \\
			\bf{Causal training (weights (\ref{eq38}))} \rule[-1ex]{0pt}{3.5ex}  & $\bf{5.08\times 10^{-4}}$ \\
			\bf{Heaviside function causal training} \rule[-1ex]{0pt}{3.5ex} & \\
			\bf{(weights (\ref{eq45}))} & $\bf{3.34 \times 10^{-4}}$ \\
			Causal training (modified MLP)~\cite{Wang2022} \rule[-1ex]{0pt}{3.5ex}  & $1.39\times 10^{-4}$ \\
			DASA-PINN~\cite{Zhang2023} \rule[-1ex]{0pt}{3.5ex}  & $8.57\times 10^{-5}$ \\
			\bf{Dirac delta function causal training} \rule[-1ex]{0pt}{3.5ex}  &  \\
			\bf{(weights (\ref{eq50}))} & $\bf{6.29\times 10^{-5}}$ \\
			\hline
		\end{tabular}
		\caption{Allen--Cahn equation: Relative ${\mathbb L}_2$ errors obtained by different approaches}
	\end{center}
\end{table}

\subsection{Korteweg--De~Vries equation}
Consider the one-dimensional Korteweg--De~Vries
equation
\begin{eqnarray}
	&& u_t + \eta u u_{x} + \mu^2 u_{xxx} = 0, \quad t \in [0,1],\quad x\in[-1,1],\label{eq64}\\
	&& u(x,0) = \cos(\pi x), \label{eq65}\\
	&& u(t,-1) = u(t,1) \label{eq66}.
\end{eqnarray}
We used the classical parameters of this problem $\eta = 1$ and $\mu = 0.022$~\cite{Zabusky1965}.
We represent the latent variable $u$ by a fully-connected neural network $u_{\pmb \theta}$ with $tanh$ activation function, 3 hidden layers and 128 neurons per hidden layer. A uniform mesh of size $100\times 512$ was constructed in computational domain $[0, 1]\times[-1, 1]$, yielding $N_{ic} = 512$ initial points and $N_r = 51200$ collocation points for enforcing the PDE residual at the algorithm~\ref{alg1} with weights given by ~(\ref{eq38}) and weights base on~\cite{Wang2022}. For the algorithm~\ref{alg2} with weights~(\ref{eq50}) the mesh is of size $100\times 256$ with $N_r = 25600$ collocation points. $\lambda_{ic} = 200$, $\lambda_r = 5$ are used in vanilla causal training~\cite{Wang2022}. For all scenarios the number of training epochs is equal to $300000$. The reference solution was generated  using the pseudo-spectral method~\cite{Zabusky1965}.

To apply algorithm~\ref{alg2} with weights~(\ref{eq50}) we need to derive the functions $\mathcal{R}\left[{\vec U}^*_{{\bm \theta}} \right]$ for current problem. This derivation is presented in the next subsection.

\subsubsection{Spatial-temporal local training}~\label{KdV}
To use the algorithm~\ref{alg2} we need to represent the value $\mathcal{R}\left[{\vec U}^*_{{\bm \theta}} \right]$ in terms of $\mathcal{R}\left[{\vec u}_{{\bm \theta}} \right]$ for $t\geq 0$. For this purpose let use expression~(\ref{eq30}).
\begin{align}
	\mathcal{R}\left[{\vec U}^*_{{\bm \theta}} \right] = \frac{\partial {\vec U}^*_{{\bm \theta}}}{\partial t} + \eta {\vec U}^*_{{\bm \theta}} \frac{\partial {\vec U}^*_{{\bm \theta}}}{\partial x} + \mu^2 \frac{\partial^3 {\vec U}^*_{{\bm \theta}}}{\partial x^3}.\label{eq67}
\end{align}

\begin{align}
	\mathcal{R}\left[{\vec U}^*_{{\bm \theta}} \right] & = \left(\frac{\partial {\vec u}_{{\bm \theta}}}{\partial t} + \eta {\vec u}_{{\bm \theta}} \frac{\partial {\vec u}_{{\bm \theta}}}{\partial x} + \mu^2 \frac{\partial^3 {\vec u}_{{\bm \theta}}}{\partial x^3}\right)\left[\chi(x - X_1) - \chi(x - X_2)\right]\notag\\
	& + \left[\frac{\partial {\vec f_1}(t)}{\partial t}\chi(X_1 - x)
	+ \frac{\partial{\vec f_2}(t)}{\partial t} \chi(x - X_2)\right]\notag\\
	& + \eta {\vec f_1}(t) \left[{\vec u}_{{\bm \theta}}(t,{x}) - {\vec f_1}(t)\right]\delta(x - X_1)\notag\\
	& + \eta {\vec f_2}(t) \left[{\vec u}_{{\bm \theta}}(t,{x}) - {\vec f_2}(t)\right]\delta(x - X_2)\notag\\
	& + 3 \mu^2 \left[\frac{\partial {\vec u}_{{\bm \theta}}}{\partial x} \frac{\partial \delta(x - X_1)}{\partial x} + \frac{\partial^2 {\vec u}_{{\bm \theta}}}{\partial x^2} \delta(x - X_1)\right]\notag\\
	& - 3 \mu^2 \left[\frac{\partial {\vec u}_{{\bm \theta}}}{\partial x} \frac{\partial \delta(x - X_2)}{\partial x} + \frac{\partial^2 {\vec u}_{{\bm \theta}}}{\partial x^2} \delta(x - X_2)\right]\notag\\
	&+ \mu^2 \left[{\vec u}_{{\bm \theta}} - {\vec f_1}(t)\right]\frac{\partial^2 \delta(x - X_1)}{\partial x^2}
	+ \mu^2 \left[{\vec f_2}(t) - {\vec u}_{{\bm \theta}}\right]\frac{\partial^2 \delta(x - X_2)}{\partial x^2}.\label{eq68}
\end{align}
By applying the fundamental equation that defines the derivatives of the delta function~\cite{Saichev1996}
\begin{equation}
	\int_{-\infty}^{\infty} f(x) \frac{\partial^{n}\delta(x)}{\partial x^n} dx = - \int_{-\infty}^{\infty} \frac{\partial f(x)}{\partial x} \frac{\partial^{n-1}\delta(x)}{\partial x^{n-1}} dx,\label{eq69}
\end{equation}
we can simplify equation~(\ref{eq68}) to form
\begin{align}
	\mathcal{R}\left[{\vec U}^*_{{\bm \theta}} \right] & = \left(\frac{\partial {\vec u}_{{\bm \theta}}}{\partial t} + \eta {\vec u}_{{\bm \theta}} \frac{\partial {\vec u}_{{\bm \theta}}}{\partial x} + \mu^2 \frac{\partial^3 {\vec u}_{{\bm \theta}}}{\partial x^3}\right)\left[\chi(x - X_1) - \chi(x - X_2)\right]\notag\\
	& + \left[\frac{\partial {\vec f_1}(t)}{\partial t}\chi(X_1 - x)
	+ \frac{\partial{\vec f_2}(t)}{\partial t} \chi(x - X_2)\right]\notag\\
	& + \eta {\vec f_1}(t) \left[{\vec u}_{{\bm \theta}}(t,{x}) - {\vec f_1}(t)\right]\delta(x - X_1)\notag\\
	& + \eta {\vec f_2}(t) \left[{\vec u}_{{\bm \theta}}(t,{x}) - {\vec f_2}(t)\right]\delta(x - X_2)\notag\\
	& + \mu^2 \frac{\partial^2 {\vec u}_{{\bm \theta}}}{\partial x^2}\left[ \delta(x - X_1) - \delta(x - X_2)\right].\label{eq70}
\end{align}
Expression (\ref{eq29}) in the domain $t\geq 0$ has the following representation
\begin{align}
	\mathcal{R}^*\left[{\vec U}_{{\bm \theta}}, {\vec g}\right] \left(t, {x} \right) &{=}
	\mathcal{R}\left[{\vec u}_{{\bm \theta}} \right] \left(t, {x} \right) {+} \left[{\vec u}_{{\bm \theta}}\left(t, {x} \right) {-} {\vec g}\left({x} \right) \right] \delta(t)\notag\\
	& + \eta {\vec f_1}(t) \left[{\vec u}_{{\bm \theta}}(t,{x}) - {\vec f_1}(t)\right]\delta(x - X_1)\notag\\
	& + \eta {\vec f_2}(t) \left[{\vec u}_{{\bm \theta}}(t,{x}) - {\vec f_2}(t)\right]\delta(x - X_2)\notag\\
	& + \mu^2 \frac{\partial^2 {\vec u}_{{\bm \theta}}}{\partial x^2}\left[ \delta(x - X_1) - \delta(x - X_2)\right]. \label{eq71}
\end{align}
For the periodic boundary conditions with minimum spatial period $X_2 - X_1$ equation~(\ref{eq71}) is simplified to
\begin{align}
	\mathcal{R}^*\left[{\vec U}_{{\bm \theta}}, {\vec g}\right] \left(t, {x} \right) &{=}
	\mathcal{R}\left[{\vec u}_{{\bm \theta}} \right] \left(t, {x} \right) {+} \left[{\vec u}_{{\bm \theta}}\left(t, {x} \right) {-} {\vec g}\left({x} \right) \right] \delta(t)\notag\\
	& + 2 \eta {\vec u}_{{\bm \theta}}(t,X_2) \left[{\vec u}_{{\bm \theta}}(t,{x}) - {\vec u}_{{\bm \theta}}(t,X_2) \right]\delta(x - X_1). \label{eq72}
\end{align}
In this case the residual loss (\ref{eq17}) is given by
\begin{align}
	\mathcal{L}_r\left(\bm \theta \right) & = \frac{1}{T}\frac{1}{V_\Omega}\int^{T}_0 dt \int_{X_1}^{X_2} d{x}  \left| \mathcal{R}\left[{\vec u}_{{\bm \theta}} \right] \left(t, {x} \right)\right|^{2} \notag\\
	& +\frac{1}{T}\frac{1}{V_\Omega}\left\{ 2 \int_{X_1}^{X_2} d{x}  \left| \mathcal{R}\left[{\vec u}_{{\bm \theta}} \right] \left(0, {x} \right) \left[{\vec u}_{{\bm \theta}}\left(0, {x} \right) {-} {\vec g}\left({x} \right) \right]\right|\right.\notag\\
	&\left. +\lim_{\substack{ \sigma \to 0 \\ t \to 0}} \frac{\sigma}{t^2 + \sigma^2}\int_{X_1}^{X_2} d{x}  \left| {\vec u}_{{\bm \theta}}\left(0, {x} \right) {-} {\vec g}\left({x} \right)\right|^{2} \right\}\notag\\
	& +\frac{1}{T}\frac{1}{V_\Omega}\left\{ 4 \eta \int_{0}^{T} d{t}  \left| \mathcal{R}\left[{\vec u}_{{\bm \theta}} \right] \left(t, X_1 \right) {\vec u}_{{\bm \theta}}(t,X_2)\left[{\vec u}_{{\bm \theta}}(t,X_1) - {\vec u}_{{\bm \theta}}(t,X_2) \right]\right|\right.\notag\\
	&+\lim_{\substack{ \sigma \to 0 \\ x \to X_1}} \frac{\sigma}{x^2 + \sigma^2} 4\eta^2\int_{0}^{T} d{t}  \left|  {\vec u}_{{\bm \theta}}(t,X_2)\left[{\vec u}_{{\bm \theta}}(t,x) - {\vec u}_{{\bm \theta}}(t,X_2) \right]\right|^{2}\notag\\
	&\left. + 4\eta^2 \lim_{\substack{ \sigma \to 0 \\ x \to X_1}} \frac{\sigma}{x^2 + \sigma^2} \lim_{\substack{ \tilde{\sigma} \to 0 \\ t \to 0}}  \frac{\tilde{\sigma}}{t^2 + \tilde{\sigma}^2}  \left|  {\vec u}_{{\bm \theta}}(t,X_2)\left[{\vec u}_{{\bm \theta}}(t,x) - {\vec u}_{{\bm \theta}}(t,X_2) \right]\right|\left| {\vec u}_{{\bm \theta}}\left(t, {x} \right) {-} {\vec g}\left({x} \right)\right| \right\}.\label{eq73}
\end{align}

\begin{figure}[t!]
	\centering
	\includegraphics[width=0.95\textwidth]{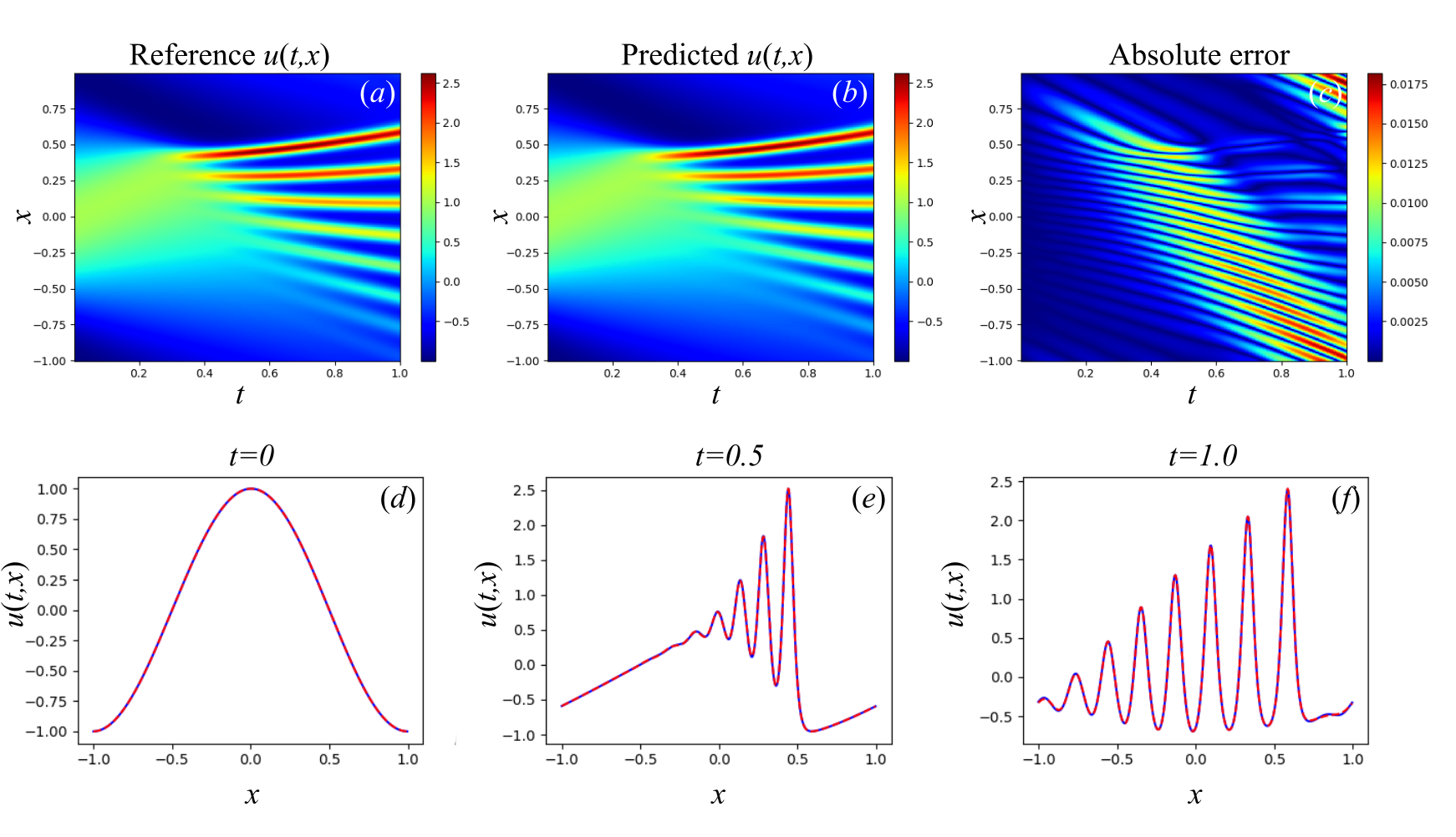}
	\caption{Korteweg--De~Vries equation. (a) is reference solution, (b) is prediction of a trained physics-informed neural network based on algorithm~\ref{alg2}, (c) is the absolute difference between reference solution and predicted one. The relative error $\epsilon_{\rm error}$ is $6.84\times 10^{-3}$. (d), (e) and (f) are comparison of the predicted (red dash lines) and reference solutions (blue solid lines) corresponding to the three temporal snapshots at $t = 0.0$, $t=0.5$ and $t=1.0$, respectively.}
	\label{fig3}
\end{figure}

\begin{table}
	\begin{center}\label{table2}
		\begin{tabular}{l|c}
			\toprule
			Method \rule[-1ex]{0pt}{3.5ex}  & Relative ${\mathbb L}_2$ error  \\
			\midrule
			Causal training (algorithm~\ref{alg1} with weights~(\ref{eq38}))  & $3.11\times 10^{-2}$ \\
			Vanilla causal training (algorithm~\ref{alg1} with weights from ~\cite{Wang2022})  & $2.03\times 10^{-2}$ \\
			{Dirac delta function causal training} & \\ (algorithm~\ref{alg2} with weights~(\ref{eq50})) & ${6.84\times 10^{-3}}$ \\
			{Dirac delta function causal training} & \\ (algorithm~\ref{alg2} with weights~(\ref{eq50}) for the PAF MLP) & ${2.45\times 10^{-3}}$ \\
			\bottomrule
		\end{tabular}
		\caption{Korteweg--DeVries equation: Relative ${\mathbb L}_2$ errors obtained by different approaches}
	\end{center}
\end{table}

Training hyperparameters are presented in the table~5 and table~6. Table~2 demonstrates the accuracy of Korteweg--De~Vries problem solution for different approaches. As can be seen the best result is achieved with Dirac delta function and algorithm~\ref{alg2} for MLP with a resulting relative ${\mathbb L}_2$ error $6.84\times 10^{-3}$. Results of this experiment are shown on the Figure~\ref{fig3}. The predicted solution coincides with the ground truth almost perfectly.

\subsubsection{MLP with heterogeneous activation functions}
Despite the relativity good results obtained for problem (\ref{eq64})--(\ref{eq66}) with the help of described fully connected neural network, the architecture discussed above is not the best for this problem. As noted in~\cite{Abbasi2022}, although the universal approximation theorem~\cite{HORNIK1989359} guarantees the existence of a neural approximator for the given problem, it does not provide any mechanism for choosing or constructing the most efficient neural network architecture for PINN training.

Let use the approach described in~\cite{Abbasi2022}, which allows to design the architecture of a neural network which demonstrates better results for this problem. This approach is named physical activation function PINN (PAF PINN) which is based on the known solutions or behavior of unknown solution $u$ for the system of differential equations under consideration.

The fundamental solution of a differential equation (\ref{eq64}) in an unbounded medium can be represented as a soliton~\cite{Korteweg1895,Vakakis2001}
\begin{align}
	u(x,t) = c\ {\rm sech}^2 \left(a x - b t - d\right),\label{eq74}
\end{align}
where $a$, $b$, $c$ and $d$ are constants, which depend on the initial condition and parameters $\eta$ and $\mu$, ${\rm sech}(\cdot) = \cosh^{-1}(\cdot) $ is hyperbolic secant. In accordance with the concept of PAF~\cite{Abbasi2022}, the activation function for one of the last layers of the neural network should be the function ${\rm sech}^2$. 

In addition, we take into account that in the system described by equation (\ref{eq64})--(\ref{eq66}) there can be several solitons whose interactions satisfy the nonlinear superposition principle~\cite{Brauer2000TheKV}. The interactions of solitons are defined by B\"{a}cklund transform which include $tanh$ function for solitons of the limited values. Therefore, as the activation function of the last layer, it is necessary to take the hyperbolic tangent.

\begin{figure}[t!]
	\centering
	\includegraphics[width=0.95\textwidth]{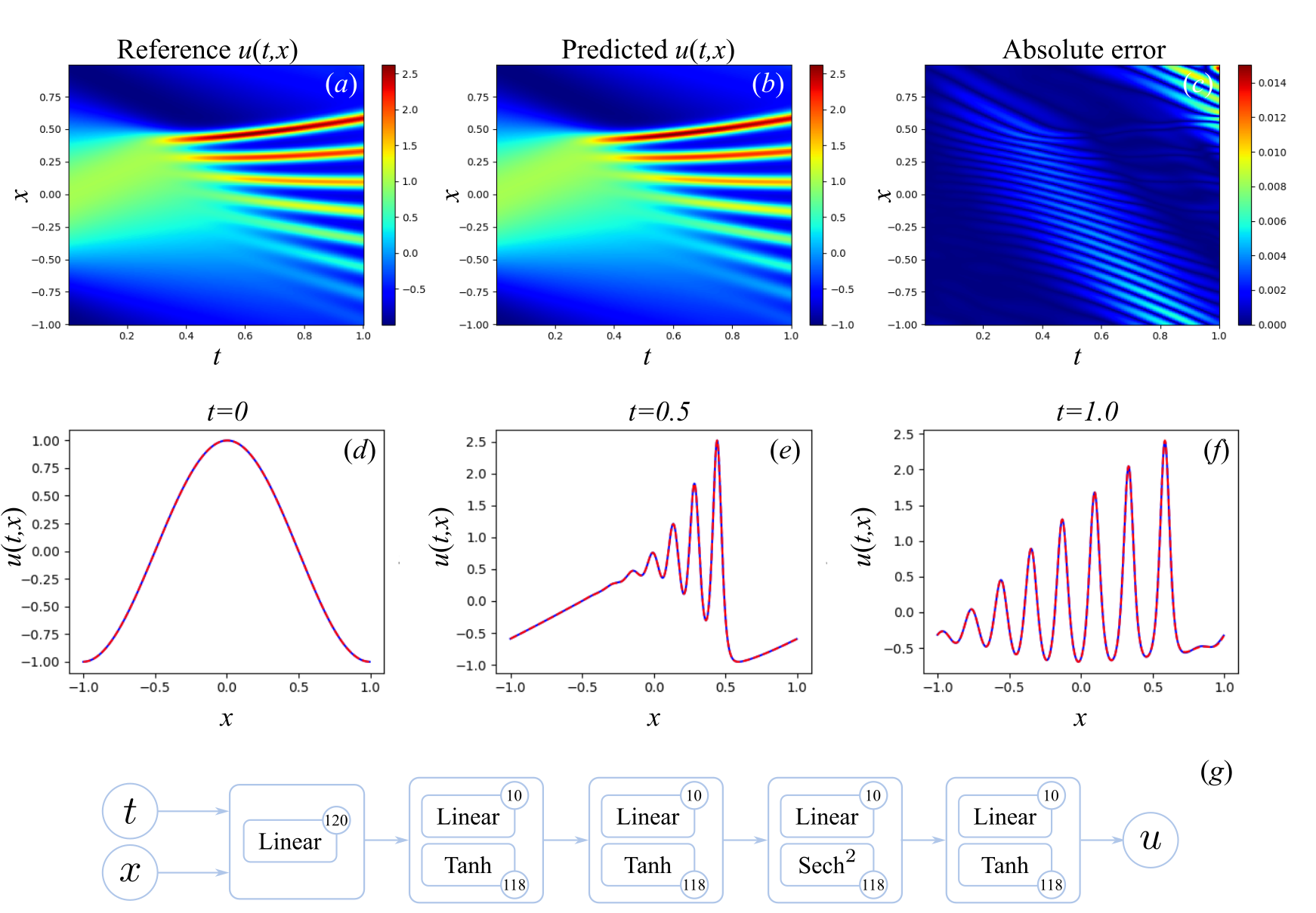}
	\caption{Korteweg--De~Vries equation. (a) is reference solution, (b) is prediction of a trained physics-informed neural network based on the algorithm~\ref{alg2}, (c) is the absolute difference between reference solution and predicted solution. The relative error $\epsilon_{\rm error}$ is $2.45\times 10^{-3}$. (d), (e) and (f) are comparison of the predicted (red dash lines) and reference solutions (blue solid lines) corresponding to the three temporal snapshots at $t = 0.0$, $t=0.5$ and $t=1.0$, respectively. (g) is schematic representation of the architecture of the neural network.}
	\label{fig4}
\end{figure}

Taking into account the noted features of the known solutions of equation (\ref{eq64}), the following neural network architecture was proposed:
\begin{enumerate}
	\item the input layer consists of 2 neurons (for $x$ and $t$ data);
	\item the second layer consists of 120 neurons with a linear activation function;
	\item the third layer consists of 128 neurons of which 10 have a linear activation function, and 118 neurons have a hyperbolic tangent as activation function;
	\item the fourth layer is the same as third;
	\item the fifth layer consists of 128 neurons of which 10 have a linear activation function, and 118 neurons have ${\rm sech}^2$ as activation function;
	\item the sixth layer is the same as third;
	\item the seventh layer consists of one linear neuron.
\end{enumerate}
A schematic representation of this neural network is shown in the panel (g) of the figure~\ref{fig4}.

Hyperparameters of Dirac delta function causal training for this neural network are presented in the table~5 and table~6. The training consists of $3$ iterations with $3\times 10^4$ epochs for each of first two iteration, and  $9\times 10^4$ epochs for the last one. The relative error $\epsilon_{\rm error}$ after training at the first, second and third iterations equals to $8.45\times 10^{-3}$, $4.03\times 10^{-3}$ and $2.45\times 10^{-3}$, respectively. Results of this experiment are shown on the Figure~\ref{fig4}. As can be seen the predicted solution coincides with the ground truth with great accuracy.

\begin{figure}[t!]
	\centering
	\includegraphics[width=0.95\textwidth]{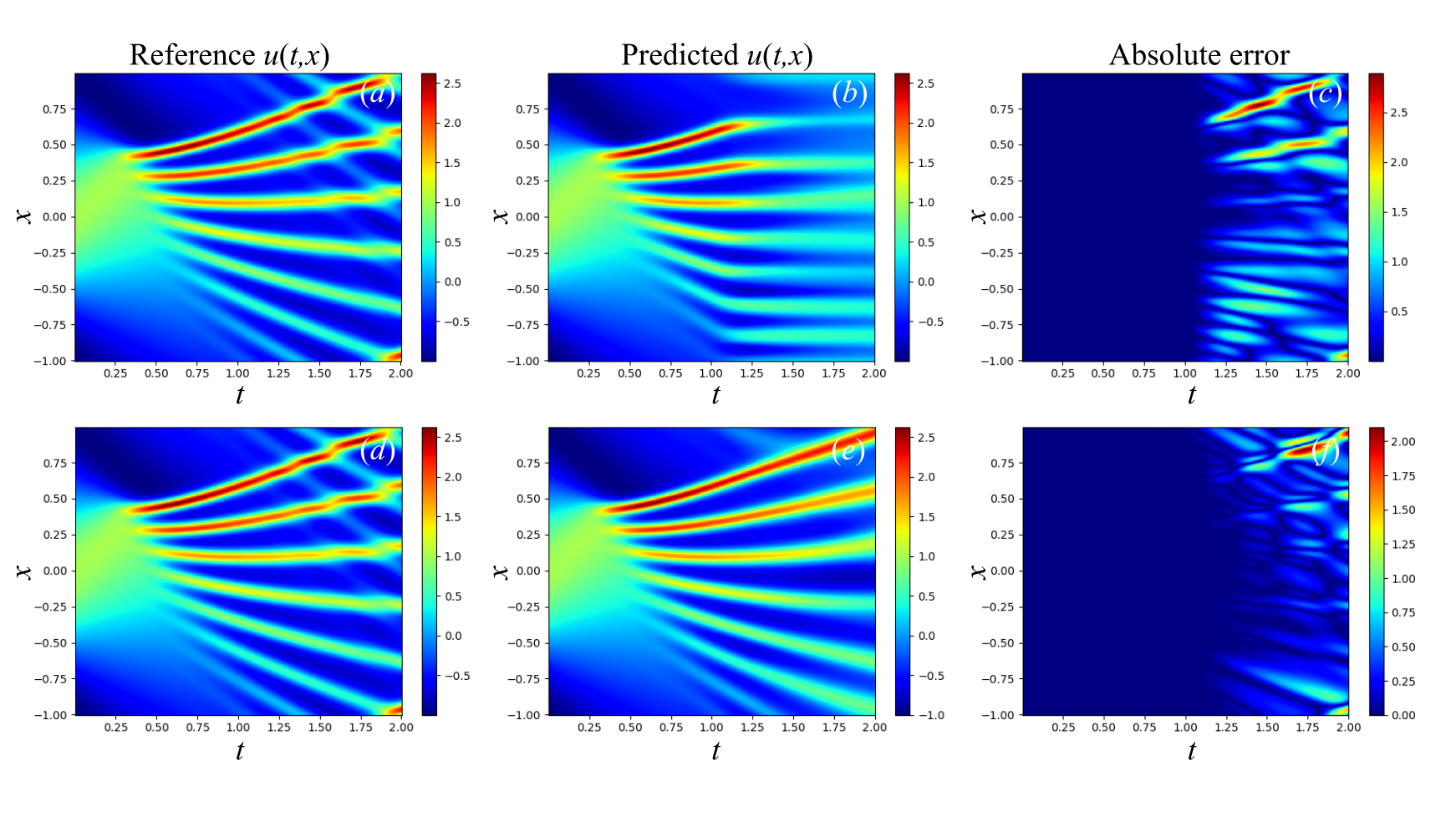}
	\caption{Korteweg--De~Vries equation. (a) is reference solution, (b) is prediction of a trained physics-informed neural	network based on algorithm~\ref{alg1} (best result for MLP, $\epsilon_{\rm error} = 6.84\times 10^{-3}$ for $t \in [0,1]$ and  $\epsilon_{\rm error} = 6.61\times 10^{-1}$ for $t \in [0,2]$), (c) is absolute value of difference of reference solution and predicted solution. (d)--(f) are same for the PAF MLP ($\epsilon_{\rm error}=2.45\times 10^{-3}$ for $t \in [0,1]$ and  $\epsilon_{\rm error} = 3.28\times 10^{-1}$ for $t \in [0,2]$).}
	\label{fig5}
\end{figure}

Based on the experiments we have come to the conclusion that the neural network with heterogeneous activation functions not only has better accuracy and speed of convergence compared with a neural network with homogeneous activation functions, but also superiorly extrapolate the solution in space-time domain where the training was not performed. Figure~\ref{fig5} shows the solutions for the MLP (b) described in the previous subsection and for PAF MLP (e) for time interval $t \in [0,2]$. It can be seen, PAF MLP demonstrates significantly better qualitative behavior of the solution in the domain $t \in [1,2]$, in comparison with the behavior of MLP in the same area. The relative error quantitatively represents these results: PAF MLP has relative error $\epsilon_{\rm error} = 3.28\times 10^{-1}$, whereas MLP  has relative error $\epsilon_{\rm error} = 6.61\times 10^{-1}$.

\subsection{Petrov--Kudrin equations}
Consider the Petrov--Kudrin system of equations~\cite{Petrov2010}, which describes the propagation of intense electromagnetic waves in a nonlinear nondispersive medium, 
\begin{eqnarray}
	&& E_\rho - {\varepsilon_1^{-1/2}} H_\tau = 0,\notag\\
	&& H_\rho + {\rho^{-1}} H - \varepsilon_1^{1/2} \exp\left(\alpha E\right) E_\tau = 0, \label{eq75}
\end{eqnarray}
where $\rho \in [0,5]$, $\tau \in [0, 4.75]$, $\alpha=1$, $\varepsilon_1=2$. We will consider the initial problem with initial values
\begin{eqnarray}
	&& E(0,\rho) = \left[1 + \rho^2 \exp\left(\alpha E\right)\right]^{-3/2},\quad H(0,\rho) = 0. \label{eq76}
\end{eqnarray}

This problem has the exact solution~\cite{Petrov2010} in the following form
\begin{eqnarray}
	&&\hspace{-1cm} E(\tau,\rho) = {\rm Re} \left\{\left[\left(1 - i \theta\right)^2 + \rho^2 \exp\left(\alpha E\right)\right]^{-3/2}\left(1 - i \theta\right)\right\},\notag\\
	&&\hspace{-1cm} H(\tau,\rho) = \varepsilon_1^{1/2}\rho \exp\left(\alpha E\right) {\rm Re} \left\{i\left[\left(1 - i \theta\right)^2 + \rho^2 \exp\left(\alpha E\right)\right]^{-3/2}\right\}. \label{eq77}
\end{eqnarray}
Here $\theta=\tau + \alpha \rho H / (2 \varepsilon_1^{1/2})$.

We represent the latent variables $E$ and $H$ with the help of a fully-connected neural network $u_{\pmb \theta}$ with $tanh$ activation function, 6 hidden layers and 64 neurons per hidden layer. For simplicity, a uniform mesh of size $100\times 256$ was constructed in the computational domain $[0, 4.75]\times[0, 5]$, yielding $N_{ic} = 512$ initial points and $N_r = 25600$ collocation points for enforcing the PDE residual in algorithm~\ref{alg1} with weights~(\ref{eq38}), (\ref{eq50}) or with weights from~\cite{Wang2022}. For all scenarios the number of epochs is supposed to be $300000$. As in the case of the Allen--Cahn problem, for this problem, the logarithm of (\ref{eq16}) and the optimization problem~
\ref{eq63} were used for training $u_{\pmb \theta}$.

\begin{table}
	\begin{center}\label{table3}
		\begin{tabular}{l|c}
			\toprule
			Method \rule[-1ex]{0pt}{3.5ex}  & Relative ${\mathbb L}_2$ error  \\
			\midrule
			Vanilla causal training (algorithm~\ref{alg1} with weights from~\cite{Wang2022})  & $8.01\times 10^{-3}$ \\
			Causal training (algorithm~\ref{alg1} with weights~(\ref{eq38}))  & $7.82\times 10^{-3}$ \\
			{Dirac delta function causal training} & \\ (algorithm~\ref{alg1} with weights~(\ref{eq50})) & ${7.69\times 10^{-3}}$ \\
			\bottomrule
		\end{tabular}
		\caption{Petrov--Kudrin equations: Relative ${\mathbb L}_2$ errors obtained by different approaches}
	\end{center}
\end{table}

To use the algorithm~\ref{alg1} it is necessary to get away from the limit representation of delta function towards the discrete form. The required transformation is described in the following subsection.

\subsubsection{Approximation of Dirac delta function }\label{approxDelta}
According to the equation~(\ref{eq22}) let approximate the delta function in the (\ref{eq23}) with function $1/\sigma$. It is assumed that the $\sigma$ value is the doubled time during which the initial distribution of the field decreases or increases twice. For the problem (\ref{eq75}),(\ref{eq76}) we need to find dependence of field $E$ on $\tau$ at $\rho = 0$. This value is further denoted with $\tilde{E}$. In addition, we suppose that the field $E$ depends on the coordinate $\rho$ similar to $E$ in (\ref{eq76}). Taking into account the latter circumstance, we obtain the following equation for the $\tilde{E}$ from the equations~(\ref{eq75})
\begin{equation}
	\frac{\partial}{\partial \tau}\left(e^{\alpha \tilde{E}} \frac{\partial \tilde{E}}{\partial \tau}\right) = -6 e^{\alpha \tilde{E}}.\label{eq78}
\end{equation}
The solution of the equation (\ref{eq78}) is
\begin{equation}\label{eq79}
	\tilde{E}(\tau) = \frac{1}{\alpha} \log\left[e^{\alpha}\cos\left(\tau\sqrt{6\alpha} \right)\right].
\end{equation}
Here we used initial value $\tilde{E} = 1$. It follows from the equation~(\ref{eq79}) that the decay time of the field $\tilde{E}$ is two times
\begin{equation}\label{eq80}
	\tilde{\tau} = \arccos\left(e^{-\alpha/2}\right) / \sqrt{6\alpha}.
\end{equation}
It follows from this equation that
\begin{equation}\label{eq81}
	\delta(\tau) \simeq 1/ (2 \tilde{\tau}).
\end{equation}

\begin{figure}[t!]
	\centering
	\includegraphics[width=0.95\textwidth]{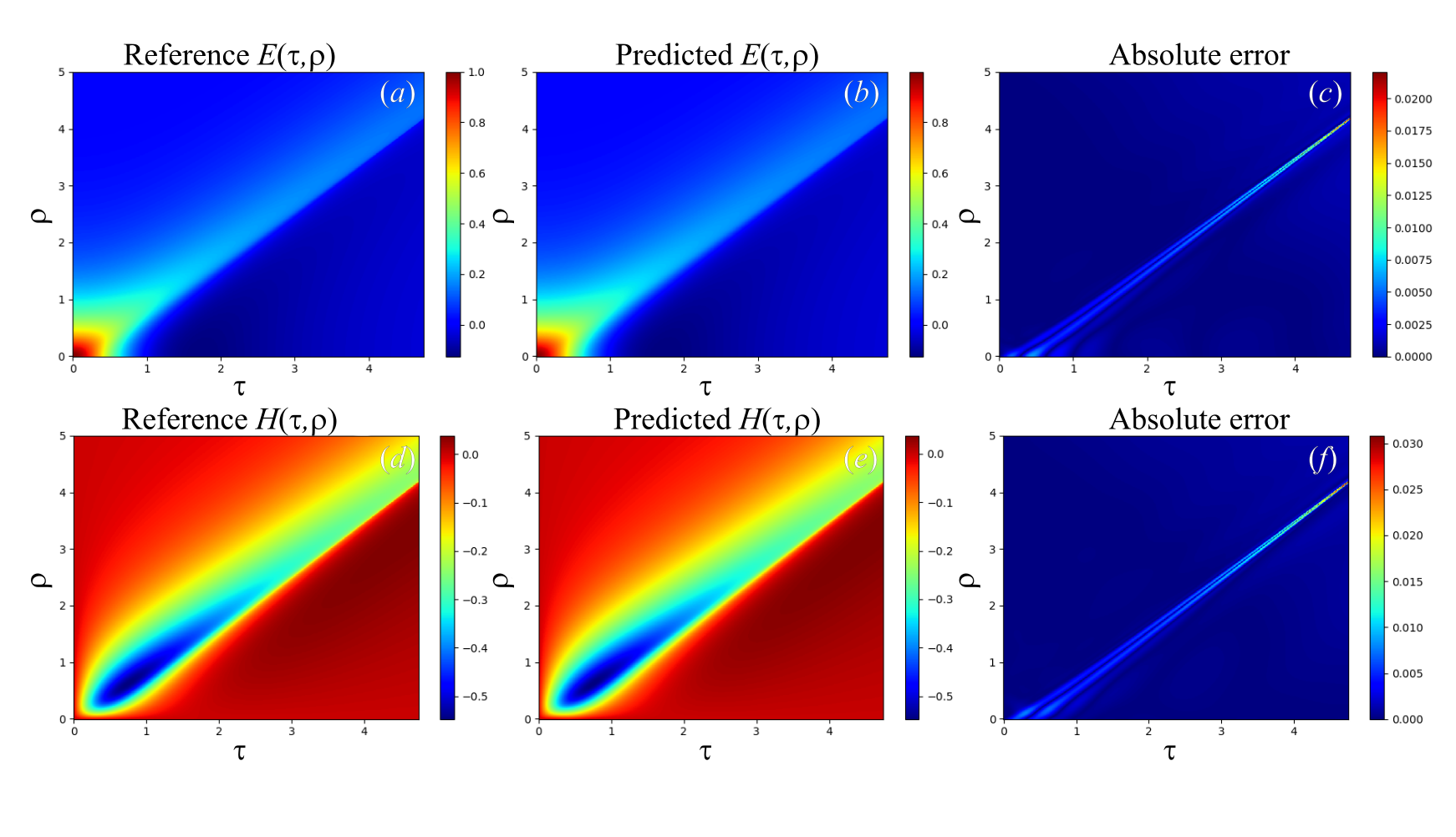}
	\caption{Petrov--Kudrin equations: the first row is $E$, the second row is $H$. (a) and (d) are reference solutions, (b) and (e) are predictions of a trained physics-informed neural network based on the algorithm~\ref{alg1}, (c) and (f) are absolute values of difference between reference and predicted solutions. The total relative error $\epsilon_{\rm error}$ is $7.69\times 10^{-3}$.}
	\label{fig6}
\end{figure}

\begin{figure}[t!]
	\centering
	\includegraphics[width=0.95\textwidth]{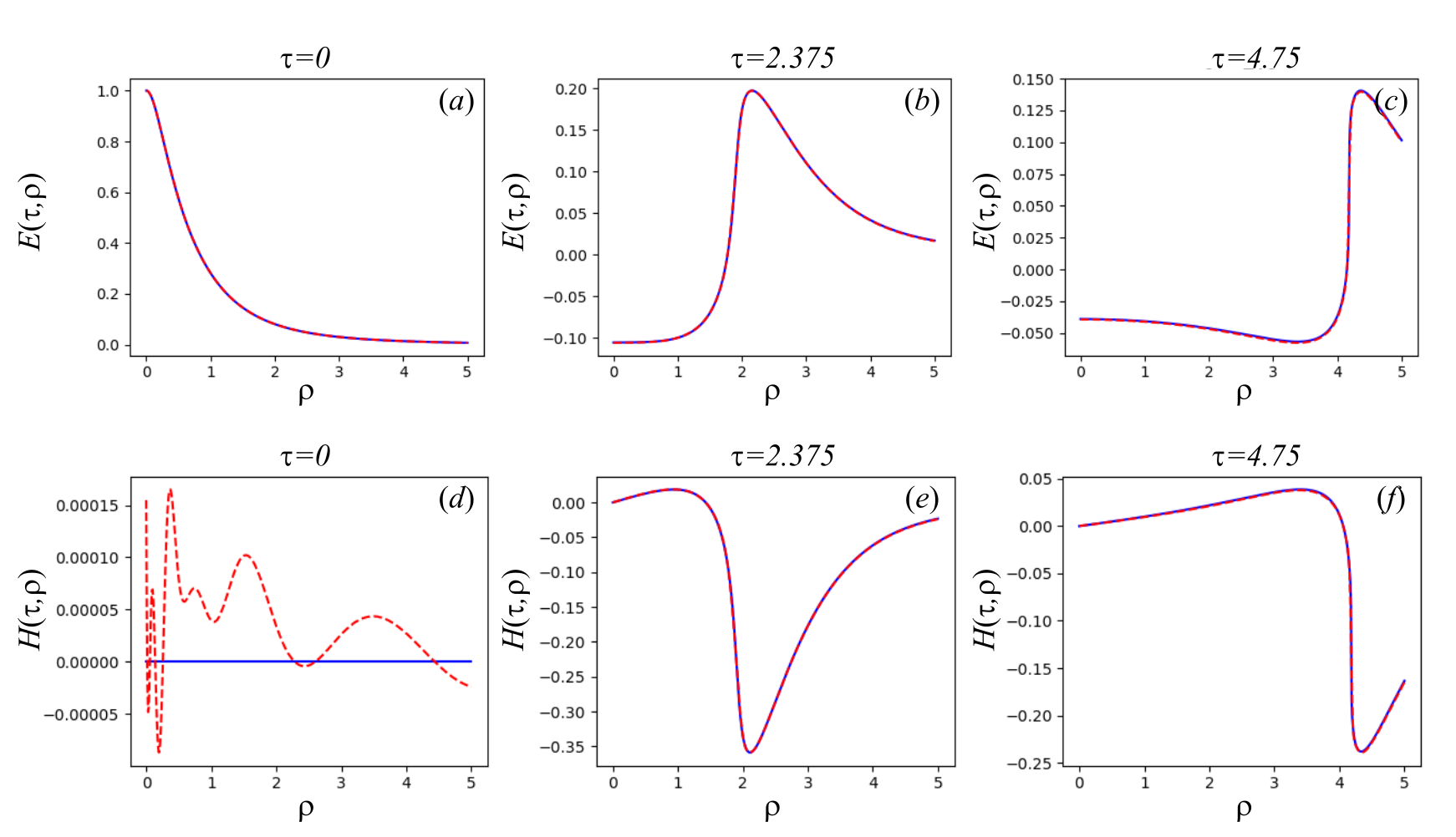}
	\caption{Petrov--Kudrin equations: (a), (b) and (c) are comparison of the predicted (red dash lines) and reference solutions (blue solid lines) of electric field $E$ corresponding to the three temporal snapshots at $\tau = 0.0$, $\tau=2.375$ and $\tau=4.75$, respectively.  (d), (e) and (f) are the same dependences for the magnetic field $H$}
	\label{fig7}
\end{figure}

Training hyperparameters for algorithm~\ref{alg1} with weights from~\cite{Wang2022}, with weights~(\ref{eq38}) and with weights~(\ref{eq50}) are presented in the table~5 and table~6. Table~3 presents the accuracy of the solution for Petrov--Kudrin problem obtained with the help of these approaches. The best result is achieved with Dirac delta function and algorithm~\ref{alg1} for MLP with relative ${\mathbb L}_2$ error $7.69\times 10^{-3}$. Results of this experiment are shown on the Figures~\ref{fig6} and \ref{fig7}. The predicted solution coincides with the ground truth with great accuracy.

	\section{Conclusion}
	
	In this article, we presented several approaches to neural network training within the framework of the PINN concept. These approaches simplify the construction and analysis of the loss function and improve the training convergence. We have also proposed the neural network architectures which are more relevant for the problems under consideration.
	
	We have reformulated the original problem described by the differential equation and the initial and boundary conditions into the problem described only by the differential equation. The major advantage of such a trick is that it becomes possible to reduce the loss function to the single term associated with the differential equations, thus eliminating the need to tune the scaling coefficients for the terms associated with the boundary and initial conditions. Based on this approach we have derived closed-form expressions for the MAE and MSE losses. Approximations of the Dirac delta function, which is part of the losses in the MSE formulation, have been proposed. These approximations are determined by the discretization method of the domain or by the behavior of the expected solution at the initial time and at the boundaries of the domain.
	
	Inspired by~\cite{Wang2022} we have proposed loss functions based on causality and spatial-temporal locality principles within the framework of the formalism of generalized functions such as Heaviside step and Dirac delta functions. It has also been shown that the causality training can be extended to the spatial domains and can be transformed into the spatial-temporal training.

	Finally, a modification of the MLP model for the Korteweg--De~Vries problem has been proposed based on~\cite{Abbasi2022}. This modification is based on the fundamental solution of the Korteweg--De~Vries equation and the nonlinear superposition principle for the solitons. It has been shown that such a neural network with heterogeneous activation functions can be trained much faster and the obtained results have significantly better accuracy compared to MLP with homogeneous activation functions. In addition, such neural network is much better able to extrapolate to new domains.
	
	Numerical experiments have been performed for a number of problems and the accuracies of the proposed methods are given. The causal trainings with weights based on the Dirac delta function were achieved the best results in all carried out the numerical experiments.

\newpage

\appendix

\section{Nomenclature}\label{appA}

The summarizes the main notations, abbreviations and symbols are given in table~4.
\begin{table}[h!]
	\centering
	\begin{tabular}{ll}
		\toprule
		Notation & Description  \\
		\midrule
		MLP & Multilayer perceptron  \\
		PDE & Partial differential equation\\
		PINN & Physic-informed neural network \\
		MAE & Mean absolute error\\
		MSE & Mean squared error \\
		DASA-PINN & Differentiable adversarial self-adaptive pointwise weighing scheme for PINN\\
		PAF & physical activation function\\
		${\vec u}$ & solution of PDE	\\
		$\bf U$ & artificially extended solution of PDE\\
		${\mathcal{N}}$ & nonlinear differential operator\\
		${\mathcal{B}}$ & boundary operator \\
		$\mathcal{R}$ & PDE residual \\
		${\vec u}_{\pmb \theta}$ & neural network representation of the PDE solution \\
		${\pmb \theta}$ & vector of the trainable parameters of the neural network \\
		$N_t$ & number of intervals mesh along the time axis\\
		$N_x$ & number of intervals mesh along the $x$ axis\\
		$w_i$ & residual weight at the time $t_i$\\
		$w_n^{(\to)}$ & residual weight at the coordinate $x_n$ in case of the beginning calculation of weight at low boundary\\
		$w_n^{(\leftarrow)}$ & residual weight at the coordinate $x_n$ in case of the beginning calculation of weight at upper boundary\\
		$\varepsilon$ & causality parameter\\
		$\delta_w$ & threshold for increasing a causality parameter\\
		${\mathcal{L}\left(t,{\pmb \theta}\right)}$ & temporal residual loss\\
		${\mathcal{L}^*\left(t,{\pmb \theta}\right)}$ & temporal residual loss in the spatial-temporal domain\\
		${\mathcal{L}(\pmb \theta)}$ & aggregate training loss\\
		\bottomrule
	\end{tabular}\label{tableAppA}
	\caption{Nomenclature}
\end{table}

\newpage
\section{Hyperparameters of approaches}

Tables~5 and 6 summarize the training and network hyperparameters, respectively, for all numerical experiments. 

The optimizer used in all experiments is Adam, and parameters $\theta$ are initialized with the Xavier scheme~\cite{Glorot10a}.

\begin{table}[h!]
	\centering
	\begin{tabular}{llllllll}
		\toprule
		PDE Problem \rule[-1ex]{0pt}{3.5ex}  & Weights & $N$& $\eta_s$ & $\eta_{\min}$ & $\varepsilon_{\rm int}$ & $\delta_w$ & Scheduler  \\
		\midrule
		\multirow{3}{*}{Allen--Cahn} & (\ref{eq38}) & \multirow{2}{*}{$3\times 10^5$} & \multirow{2}{*}{$3\times 10^{-3}$} & $10^{-12}$ & \multirow{2}{*}{$10^{-3}$} & $0.99$ & \multirow{2}{*}{ExponentialLR} \\	
		& (\ref{eq51}) & & & $10^{-5}$  & & $0.95$ & \\
		& (\ref{eq50}) & $6\times 10^5$ & $10^{-3}$ & --- & $10^4$ & $0.99$ & StepLR$^{0.95}_{5000}$ \\		
		\midrule
		\multirow{3}{*}{Korteweg--De~Vries} & \cite{Wang2022} & \multirow{3}{*}{$3\times 10^5$} & \multirow{3}{*}{$3\times 10^{-3}$} & \multirow{3}{*}{---} & \multirow{3}{*}{$10^{-2}$} & \multirow{3}{*}{$0.99$} & \multirow{3}{*}{StepLR$^{0.9}_{5000}$} \\
		& (\ref{eq38}) & & & & & \\
		& (\ref{eq50}) & & & & & \\
		\midrule
		\multirow{3}{*}{Korteweg--De~Vries (PAF)} & \multirow{3}{*}{(\ref{eq50})} & \multirow{2}{*}{$3\times 10^4$} & {$3\times 10^{-3}$} & \multirow{3}{*}{$10^{-8}$} & {$10^{-5}$} & \multirow{3}{*}{$0.99$} & \multirow{3}{*}{CosineAnnealingLR} \\
		&  & & $10^{-4}$& & $0.16$ & \\
		& & $9\times 10^4$ & $3\times 10^{-5}$& & $0.33$ & \\
		\midrule
		\multirow{3}{*}{Petrov--Kudrin} & \cite{Wang2022} & \multirow{3}{*}{$3\times 10^5$} & \multirow{3}{*}{$10^{-3}$} & \multirow{3}{*}{---} & \multirow{3}{*}{$10^{-2}$} & \multirow{3}{*}{$0.99$} & \multirow{3}{*}{StepLR$^{0.9}_{5000}$} \\
		& (\ref{eq38}) & & & & & \\
		& (\ref{eq50}) & & & & & \\
		\bottomrule
	\end{tabular}\label{tableAppC1}
	\caption{Training hyperparameters}
\end{table}
Here $N$ is number of epochs, $\eta_s$ is initial value of learning rate,  $\eta_{\min}$ is minimum value of learning rate, $\varepsilon_{\rm int}$ is the initial value of the causality parameter. For the ExponentialLR scheduler the a decay--rate is calculated with following formula
\begin{equation}
	\gamma= \left(\eta_{\min}/\eta_s\right)^{1/N}.\notag
\end{equation}
StepLR$^{\gamma}_{E}$ means decay--rate is $\gamma$ every $E$ training epochs.

\begin{table}[h!]
	\centering
	\begin{tabular}{llllllll}
		\toprule
		PDE Problem \rule[-1ex]{0pt}{3.5ex}  & Architecture & \head{1.5cm}{Hidden layers}& \head{1.5cm}{Neurons per layer} & \head{1.5cm}{Activation function} & $N_t$ & $N_x$  \\
		\midrule
		{Allen--Cahn} & MLP & 4 & 128 & $\tanh$ & 100 & 256 \\		
		\midrule
		{Korteweg--De~Vries} & MLP & 3 & 128 & $\tanh$ & 100 & 512 for \cite{Wang2022} and (\ref{eq38})  \\
		& & & & & & 256 for (\ref{eq50}) \\
		\midrule
		{Korteweg--De~Vries} & PAF MLP & 5 & 120 or 128 & linear, $\tanh$ or ${\rm sech}^2$ & 100 & 512 for (\ref{eq50}) \\
		\midrule
		Petrov--Kudrin & MLP & 6 & 64 & $\tanh$ & 100 & 256 \\		
		\bottomrule
	\end{tabular}\label{tableAppC2}
	\caption{The models and mesh parameters}
\end{table}
Here $N_t$ and $N_x$ are number of values along time and $x$ coordinate axes, respectively.

\newpage
\bibliographystyle{IEEEtran}
\bibliography{Eskin_causal_loss}

\end{document}